\documentclass[twoside, a4paper, 10pt]{amsart}
\title{Multidimensional lower density versions of Pl\"unnecke's inequality}
\author{Kamil Bulinski}

\usepackage{amsfonts}
\usepackage{amsthm}
\usepackage{verbatim}
\usepackage{amsmath, amssymb}
\usepackage{tikz}
\usetikzlibrary{matrix, arrows}
\usepackage{listings}
\usepackage{color}
\usepackage{listings}
\usepackage[all]{xy}
\usepackage[pdftex,colorlinks,linkcolor=blue,citecolor=blue]{hyperref}
\usepackage{graphicx}
\usepackage{float}
\usepackage[margin=3cm]{geometry}
\usepackage{bigints}
\usepackage[toc,page]{appendix}
\usepackage{enumitem}
\allowdisplaybreaks
\begin{document}
\maketitle
\raggedbottom

%% Mathcal large
\newcommand{\cA}{\mathcal{A}}
\newcommand{\cB}{\mathcal{B}}
\newcommand{\cC}{\mathcal{C}}
\newcommand{\cD}{\mathcal{D}}
\newcommand{\cE}{\mathcal{E}}
\newcommand{\cF}{\mathcal{F}}
\newcommand{\cG}{\mathcal{G}}
\newcommand{\cH}{\mathcal{H}}
\newcommand{\cI}{\mathcal{I}}
\newcommand{\cJ}{\mathcal{J}}
\newcommand{\cK}{\mathcal{K}}
\newcommand{\cL}{\mathcal{L}}
\newcommand{\cM}{\mathcal{M}}
\newcommand{\cN}{\mathcal{N}}
\newcommand{\cO}{\mathcal{O}}
\newcommand{\cP}{\mathcal{P}}
\newcommand{\cQ}{\mathcal{Q}}
\newcommand{\cR}{\mathcal{R}}
\newcommand{\cS}{\mathcal{S}}
\newcommand{\cT}{\mathcal{T}}
\newcommand{\cU}{\mathcal{U}}
\newcommand{\cV}{\mathcal{V}}
\newcommand{\cW}{\mathcal{W}}
\newcommand{\cX}{\mathcal{X}}
\newcommand{\cY}{\mathcal{Y}}
\newcommand{\cZ}{\mathcal{Z}}
%% Mathbb large
\newcommand{\bA}{\mathbb{A}}
\newcommand{\bB}{\mathbb{B}}
\newcommand{\bC}{\mathbb{C}}
\newcommand{\bD}{\mathbb{D}}
\newcommand{\bE}{\mathbb{E}}
\newcommand{\bF}{\mathbb{F}}
\newcommand{\bG}{\mathbb{G}}
\newcommand{\bH}{\mathbb{H}}
\newcommand{\bI}{\mathbb{I}}
\newcommand{\bJ}{\mathbb{J}}
\newcommand{\bK}{\mathbb{K}}
\newcommand{\bL}{\mathbb{L}}
\newcommand{\bM}{\mathbb{M}}
\newcommand{\bN}{\mathbb{N}}
\newcommand{\bO}{\mathbb{O}}
\newcommand{\bP}{\mathbb{P}}
\newcommand{\bQ}{\mathbb{Q}}
\newcommand{\bR}{\mathbb{R}}
\newcommand{\bS}{\mathbb{S}}
\newcommand{\bT}{\mathbb{T}}
\newcommand{\bU}{\mathbb{U}}
\newcommand{\bV}{\mathbb{V}}
\newcommand{\bW}{\mathbb{W}}
\newcommand{\bX}{\mathbb{X}}
\newcommand{\bY}{\mathbb{Y}}
\newcommand{\bZ}{\mathbb{Z}}

\newcounter{dummy} \numberwithin{dummy}{section} 

\numberwithin{equation}{section}

\theoremstyle{definition}
\newtheorem{mydef}[dummy]{Definition}
\newtheorem{prop}[dummy]{Proposition}
\newtheorem{corol}[dummy]{Corollary}
\newtheorem{thm}[dummy]{Theorem}
\newtheorem{lemma}[dummy]{Lemma}
\newtheorem{eg}[dummy]{Example}
\newtheorem{notation}[dummy]{Notation}
\newtheorem{remark}[dummy]{Remark}
\newtheorem{claim}[dummy]{Claim}
\newtheorem{Exercise}[dummy]{Exercise}
\newtheorem{question}[dummy]{Question}

\begin{abstract}We investigate the lower asymptotic density of sumsets in $\bN^2$ by proving certain Pl\"unnecke type inequalities for various notions of lower density in $\bN^2$. More specifically, we introduce a notion of lower \textit{tableaux} density in $\bN^2$ which involves averaging over convex tableaux-shaped regions in $\bN^2$ which contain the origin.  This generalizes the well known Pl\"unnecke type inequality for the lower asymptotic density of sumsets in $\bN$. We also provide a conjectural Pl\"unnecke inequality for the more basic notion of lower \textit{rectangular} asymtpotic density in $\bN^2$ and prove certain partial results. \end{abstract}

\tableofcontents

\section{Introduction}

\subsection{Background} 

Pl\"unnecke's classical work \cite{Plunnecke1970} provided influential techniques for studying the cardinality of \textit{sumsets} and \textit{iterated sumsets}. We recall that for subsets $A,B$ of an abelian group $(G,+)$ we can define the \textit{sumset} $$A+B=\{a+b \text{ }| \text{ } a \in A, b \in B \}$$ and, for positive integers $k$, the \textit{$k$-fold iterated sumset} $$kA=\{a_1 + \cdots + a_k \text{ } | \text{ }a_1, \ldots, a_k \in A\}.$$  

\begin{thm}[Pl\"unnecke \cite{Plunnecke1970}] \label{Plunnecke sumset inequality} Suppose that $A,B$ are finite subsets of some abelian group $(G,+)$ and $k\geq 1$. Then $$|A+B| \geq |A|^{1-\frac{1}{k}}|kB|^{1/k}.$$

\end{thm}

Pl\"unnecke used these techniques to improve a result of Erd\H{o}s concerning lower bounds for the Schnirelmann density of a sumset $A+B$, where $A,B \subset \bZ_{\geq 0}$ and $B$ is a basis of order $k$ (that is, $kB=\bN=\{0,1,2,\ldots \}$). We define the \textit{Schnirelmann density} of a subset $A \subset \bN=\{0,1, \ldots \}$ to be \begin{align} \label{def: sigma A} \sigma(A)=\inf_{n \in \bN} \frac{|A \cap [0,n]|}{n+1}. \end{align} Note that the usual definition is $\sigma(A)=\inf_{n \geq 1} \frac{|A \cap [1,n]|}{n}$, which is necessary for some applications such as Mann's Theorem, but the results of interest to us remain valid when using our definition. Erd\"os proved that $$\sigma(A+B)\geq \left(1+\frac{1-\sigma(A)}{2k}\right) \sigma(A)$$ for $A,B \subset \bN$ such that $kB=\bN$ for some positive integer $k$ \cite{ErdosBasis}. Pl\"unnecke greatly improved this lower bound by proving the following extension of his cardinality estimate (Theorem~\ref{Plunnecke sumset inequality}) to Schnirelmann density. 

\begin{thm}[Pl\"unnecke's inequality for Schnirelmann density \cite{Plunnecke1970}] \label{Plunnecke for Schnirelmann} For positive integers $k$ and $A,B \subset \bN$ with $0 \in B$, we have that $$\sigma(A+B) \geq \sigma(A)^{1-\frac{1}{k}}\sigma(kB)^{1/k}.$$

\end{thm}

A good account of this as well as a proof of Theorem~\ref{Plunnecke sumset inequality} and related results can be found in Ruzsa's book \cite{Ruzsasumsetsandstructure}. Although Schnirelmann density has played an important role in additive number theory (see, for instance, Schnirelmann's proof that the primes are an asymptotic basis \cite{SchnirelmannEigenschaften}), it lacks many asymptotic features such as translation invariance. From a combinatorial perspective, the \textit{lower asymptotic density}, given by $$\underline{d}(A)=\liminf_{n\to \infty} \frac{|A \cap [0,n]|}{n+1},$$ is a more natural notion of the asymptotic size of a set $A \subset \bN$. It turns out that Theorem~\ref{Plunnecke for Schnirelmann} is also true with $\underline{d}$ in place of $\sigma$. 

\begin{thm}[See \cite{JinPlunnecke}, \cite{JinEpsilon} and \cite{Ruzsasumsetsandstructure}] Suppose that $A,B \subset \mathbb{N}$ and $k\geq 1$. Then $$\underline{d}(A+B) \geq \underline{d}(A)^{1-\frac{1}{k}}\underline{d}(kB)^{1/k}.$$ \end{thm}

This was first obtained by Ruzsa \cite{Ruzsasumsetsandstructure} and an alternative proof, which the author has found insightful, was given by Jin in \cite{JinPlunnecke} and \cite{JinEpsilon}.

\subsection{Density Pl\"unnecke inequalites in (semi)groups}

It is natural to ask whether density versions of Pl\"unnecke's inequality hold in other countable abelian (semi)groups with certain notions of asymptotic density. Let us briefly mention some recently established results in this direction. We first recall a standard way of extending the notion of density to other groups that involves replacing the sequence of intervals $\left([0,N) \cap \bZ\right)_{N=1}^{\infty}$ with a sequence of asymptotically invariant finite sets.

\begin{mydef}[Densities along F{\o}lner sequences]\label{def: Folner density} Let $(G,+)$ be a countable abelian semigroup. A F{\o}lner sequence is a sequence $F_1,F_2, \ldots $ of finite subsets of $G$ that is \textit{asymptotically invariant}, i.e., for each $g \in G$ we have that $$\lim_{n \to \infty} \frac{|F_n \cap (g+F_n)|}{|F_n|} = 1.$$ Moreover, if $A \subset G$ then we define the lower asymptotic density along $(F_n)$ as $$\underline{d}_{(F_n)} (A)= \liminf_{n \to \infty} \frac{|F_n \cap A|}{|F_n|}.$$ Similairly, we may define the upper asymptotic density along $(F_n)$ as $$\overline{d}_{(F_n)} (A)= \limsup_{n \to \infty} \frac{|F_n \cap A|}{|F_n|}.$$  If $\mathcal{F}$ is a collection of F{\o}lner sequences in $G$ then we can define the lower and upper densities with respect to this collection as $$\underline{d}_{\mathcal{F}} (A) = \inf\left\{ \underline{d}_{(F_n)}(A)\text{ } | \text{ } (F_n) \in \mathcal{F} \right\}$$ and $$\overline{d}_{\mathcal{F}} (A) =\sup\left\{ \overline{d}_{(F_n)}(A)\text{ } |\text{ } (F_n) \in \mathcal{F} \right\}.$$ Finally, the lower and upper Banach densities in $G$ may be defined, respectively, as \begin{align*} d_{\ast}=\underline{d}_{\textbf{F{\o}lner}(G)} \text{ and } d^{\ast}=\overline{d}_{\textbf{F{\o}lner}(G)} \end{align*} where $\textbf{F{\o}lner}(G)$ denotes the collection of all F{\o}lner sequences in $G$. 

\end{mydef} 

\begin{thm}[\cite{BulinskiFishPlunnecke}; $k'=1$ case obtained in \cite{BjorklundFishPlunnecke}] \label{BBF Plunnecke} Suppose that $(G,+)$ is a countable abelian group and $A,B \subset G$. Then for integers $0 < k' < k$ we have 
$$d^*(A+k'B) \geq d^*(kB)^{\frac{k'}{k}} d^{*}(A)^{1-\frac{k'}{k}}$$ and $$d_*(A+k'B) \geq d^*(kB)^{\frac{k'}{k}} d_{*}(A)^{1-\frac{k'}{k}}.$$

\end{thm}

For the semigroup $G=\mathbb{N}$, the $k'=1$ cases of these inequalities were obtained by Jin in \cite{JinEpsilon}. We remark that the proofs of Theorem~\ref{BBF Plunnecke} (for arbitrary countable abelian groups $G$) use ergodic theory and it is unclear whether such techniques can be applied to densities associated to smaller classes of F{\o}lner sequences, such as lower asymptotic density. The purpose of this paper is to extend the Pl\"unnecke type inequality for lower asymptotic density to the semigroup $\mathbb{N}^2$. 

\subsection{F{\o}lner sequences and lower asymptotic density in $\bN^2$}

%We begin with some general remarks about $\textbf{F{\o}lner}(G)$, the set of all F{\o}lner sequences in a countable abelian semigroup $G$. Given a F{\o}lner sequence $F=(F_n)_{n=1}^{\infty}$ in $G$, we let $$\cS(F)= \left \{ (F_{k_n})_{n=1}^{\infty} \text{ } | \text{ } k_1<k_2< \ldots \in \bZ_{>0} \right \} \subset \textbf{F{\o}lner}(G) $$ denote the set of all subsequences of $F$. It is clear that $\underline{d}_{\cS(F)}=\underline{d}_F$. Now let $\cF \subset \textbf{F{\o}lner}(G)$ be a collection of F{\o}lner sequences closed under taking subsequences (i.e., $F \in \cF$ implies that $\cS(F) \subset \cF$). Then the collection of all \textit{rectangular} F{\o}lner sequences $\mathcal{F} \oplus \mathcal{F} = \left \{ F_n \times F'_n \text{ } | \text{ } (F_n), (F'_n) \in \mathcal{F} \right\}$ satisfies the property that the corresponding lower asymptotic density is a \textit{product density} in the sense that $$\underline{d}_{\mathcal{F}}(A) \underline{d}_{\mathcal{F}}(B)=\underline{d}_{\mathcal{F} \oplus \mathcal{F}} (A \times B)$$ for $A,B \subset G$ (note that this equality relies on $\mathcal{F}$ being closed under taking subsequences). Unlike in the case of countably additive probabilty measures, such a product density is not unique. This can be seen by taking 

A natural candidate for \textit{lower asymptotic density} in $\bN^2$ arises from considering the family $$\textbf{Rect}=\left \{ \left( [0,W_k] \times [0,H_k] \cap \mathbb{N}^2 \right)_{k=1}^{\infty} \text{ }|\text{ } W_k,H_k \to \infty \text{ as } k \to \infty \right \} \subset \textbf{F{\o}lner}(\bN^2),$$ as the corresponding lower density is a product density in the sense that $$\underline{d}_{\textbf{Rect}}(A \times B)=\underline{d}(A)\underline{d}(B),$$ for $A,B \subset \bN^2$. In the one-dimensional case, the collection $\{[0,W_k] \cap \bN^2 \text{ }| \text{ } W_k \to \infty \text{ as } k \to \infty \}$ of sequences of intervals which gives rise to the lower density $\underline{d}$ in $\bN$ satisfies the desirable property of being closed under pointwise unions\footnote{The pointwise union of two F{\o}lner sequences $(F_n)_n$ and $(G_n)_n$ is the F{\o}lner sequence $(F_n \cup G_n)_n$.}. Unfortunately, $\textbf{Rect}$ does not satisfy this property, so it seems natural to consider \begin{align}\label{def of Tab} \textbf{Tab} = \bigcup_{L=1}^{\infty} \left \{ \left(R^{(1)}_k \cup R^{(2)}_k \cup \ldots \cup R^{(L)}_k \right)_{k=1}^{\infty} \text{ } | \text{ } \left( R^{(i)}_{k}\right)_{k=1}^{\infty} \in \textbf{Rect} \text{ for } i=1,2 \ldots, L \right \}, \end{align} which is the smallest subset of $\textbf{F{\o}lner}(G)$ that contains $\textbf{Rect}$ and is closed under pointwise unions. The corresponding lower density $\underline{d}_{\textbf{Tab}}$ has the curious property that it is also a product density as above and hence agrees with $\underline{d}_{\textbf{Rect}}$ on cartesian products, i.e., $$\underline{d}_{\textbf{Tab}}(A \times B)=\underline{d}(A)\underline{d}(B)=\underline{d}_{\textbf{Rect}}(A \times B),$$ for $A,B \subset \bN$ (See Appendix~\ref{appendix: density of cartesian products}). In fact, this property is also satisfied by the collection $$\textbf{Tab}(L)=  \left \{ \left(R^{(1)}_k \cup R^{(2)}_k \cup \ldots \cup R^{(L)}_k \right)_{k=1}^{\infty} \text{ } | \text{ } \left( R^{(i)}_{k}\right)_{k=1}^{\infty} \in \textbf{Rect} \text{ for } i=1,2 \ldots, L \right \}$$ for all $L \in \bN$, which is itself a family of F{\o}lner sequences that we will be interested in.

\subsection{Lower density versions of Pl\"unnecke's inequality in $\bN^2$}

The main goal of this paper is to address the following question. 

\begin{question}\label{question: plunnecke for d_rect} For $A,B \subset \bN^2$, with $(0,0) \in B$, and positive integers $k'<k$, is it true that $$\underline{d}_{\textbf{Rect}}(A+k'B) \geq \underline{d}_{\textbf{Rect}}(A)^{1 - \frac{k'}{k}}\underline{d}_{\textbf{Rect}}(kB)^{\frac{k'}{k}}?$$

\end{question}

Our first partial result is an affirmative answer to this question when $A \subset \mathbb{N}^2$ is such that the density $d_{\textbf{Rect}}(A)$ \textit{exists}, by which we mean that $$\underline{d}_{\textbf{Rect}}(A)=\overline{d}_{\textbf{Rect}}(A).$$ In this case, we denote this common value by $d_{\textbf{Rect}}(A)$.

\begin{prop}[Pl\"unnecke inequalities for $\underline{d}_{\textbf{Rect}}$ when density of $A$ exists] \label{prop: plunnecke when d exists} Suppose that $A \subset \mathbb{N}^2$ is such that $d_{\textbf{Rect}}(A)$ exists. Then for all $B \subset \mathbb{N}^2$ with $(0,0) \in B$ and integers $0<k'<k$ we have that $$\underline{d}_{\textbf{Rect}}(A+k'B) \geq d_{\textbf{Rect}}(A)^{1 - \frac{k'}{k}}\underline{d}_{\textbf{Rect}}(kB)^{\frac{k'}{k}}.$$ 

\end{prop}

\begin{remark} Note that it is very easy to affirmatively answer Question~\ref{question: plunnecke for d_rect} up to a constant by using the fact that $(A \cap [0,N] \times [0,M]) + (B \cap [0,N] \times [0,M]) \subset (A+B) \cap [0,2N] \times [0,2M]$. For instance, one may use this fact to immediately deduce that for $A,B \subset \mathbb{N}^2$ and positive integers $k$, we have that $$\underline{d}_{\textbf{Rect}}(A+B) \geq \frac{1}{4} \underline{d}_{\textbf{Rect}}(A)^{1 - \frac{1}{k}}\underline{d}_{\textbf{Rect}}(kB)^{\frac{1}{k}}.$$

The constant becomes much worse than $\frac{1}{4}$ if one applies this method to estimate the lower densities of $A+k'B$ for $k'<k$ large; more precisely, one gets $$\underline{d}_{\textbf{Rect}}(A+k'B) \geq \frac{1}{(1+k')^2} \underline{d}_{\textbf{Rect}}(A)^{1 - \frac{k'}{k}}\underline{d}_{\textbf{Rect}}(kB)^{\frac{k'}{k}}.$$

\end{remark}

\subsection{Statements of main results and applications} 

We are able to affirmatively answer Question~\ref{question: plunnecke for d_rect} if we replace $\textbf{Rect}$ with the collection $\textbf{Tab}$ introduced above in (\ref{def of Tab}). An element of $\textbf{Tab}$ will be refered to as a \textit{tableau F{\o}lner sequence} and we will refer to the corresponding notion of lower density $\underline{d}_{\textbf{Tab}}$ (as per Definition~\ref{def: Folner density}) as the \textit{lower tableau density}.

\begin{thm}\label{thm: plunnecke for d_tab} Let $A,B \subset \mathbb{N}^2$ such that $(0,0) \in B$. Then $$\underline{d}_{\textbf{Tab}}(A+k'B) \geq \underline{d}_{\textbf{Tab}}(A)^{1-k'/k} \underline{d}_{\textbf{Tab}}(kB)^{k'/k}$$ for integers $0<k'<k$.

\end{thm} 

In fact, our techniques also give the following partial answer to Question~\ref{question: plunnecke for d_rect}.

\begin{thm}\label{thm: drect bound in terms of dtab} Suppose that $0<k'<k$ are integers and $A, B \subset \mathbb{N}^2$ such that $(0,0) \in B$, then $$ \underline{d}_{\textbf{Rect}}(A+k'B) \geq \underline{d}_{\textbf{Tab}}(A)^{1 - \frac{k'}{k}}\underline{d}_{\textbf{Rect}}(kB)^{\frac{k'}{k}}.$$

\end{thm}

This enables us to affirmatively answer Question~\ref{question: plunnecke for d_rect} for a broader class of examples not covered by Proposition~\ref{prop: plunnecke when d exists}.

\begin{corol}\label{corol: d_rect is d_tab} Suppose that $0<k'<k$ are integers and $A, B \subset \mathbb{N}^2$ such that $(0,0) \in B$, $\underline{d}_{\textbf{Tab}}(A) = \underline{d}_{\textbf{Rect}}(A)$. Then $$\underline{d}_{\textbf{Rect}}(A+k'B) \geq \underline{d}_{\textbf{Rect}}(A)^{1-\frac{k'}{k}}\underline{d}_{\textbf{Rect}}(kB)^{\frac{k'}{k}}.$$

\end{corol}

In particular, since $\underline{d}_{\textbf{Tab}}(A)=d_{\textbf{Rect}}(A)$ whenever $d_{\textbf{Rect}(A)}$ exists, we get Proposition~\ref{prop: plunnecke when d exists}. Theorems~\ref{thm: plunnecke for d_tab} and \ref{thm: drect bound in terms of dtab} are both immediate consequences of the following general result. 

\begin{thm} \label{main result in intro} Suppose that $0<k'<k$ are integers and $A, B \subset \mathbb{N}^2$ such that $(0,0) \in B$, then $$ \underline{d}_{\textbf{Tab}(L)}(A+k'B) \geq \underline{d}_{\textbf{Tab}}(A)^{1 - \frac{k'}{k}}\underline{d}_{\textbf{Tab}(L)}(kB)^{\frac{k'}{k}}$$ for all $L \in \bN$. \end{thm}

\subsection{Examples: Fractal Sets} We now turn to constructing some examples of subsets of $\mathbb{N}^2$ which demonstrate the novelty of our main result and its corollary. We are able to give lower bounds for $\underline{d}_{\textbf{Rect}}(A+k'B)$ in the case where $A$ possesses a certain fractal structure and $B$ is a rectangular asymptotic basis of order $k>k'$ (that is, $\underline{d}_{\textbf{Rect}}(kB) = 1$). We give an example of such a \textit{fractal set} before giving a broad definition.

\begin{eg}[A fractal set and an application of Theorem~\ref{thm: plunnecke for d_tab}]\label{eg: fractal set example} In this example, by $[a,b)$ we mean $\{x \in \mathbb{N} | a \leq x < b\}$ for $a<b \in \mathbb{N}$. Choose a sequence $0<u_1<u_2< \ldots$ in $\mathbb{N}^2$ such that $u_{k}>2u_{k-1}$ for $k>1$ and $\lim_{k \to \infty} \frac{u_k}{u_{k-1}} = \infty$. Let $A \subset \mathbb{N}$ be the set given in Figure~\ref{Fig: fractal example}, more precisely $$A=[0,u_1)^2 \cup [u_1,2u_1)^2 \cup \bigcup_{k=2}^{\infty} [u_k,2u_k)^2 \cup \left( [0,u_{k})^2 \setminus [0,2u_{k-1})^2\right). $$

\begin{figure}[H]
\includegraphics[scale=0.35]{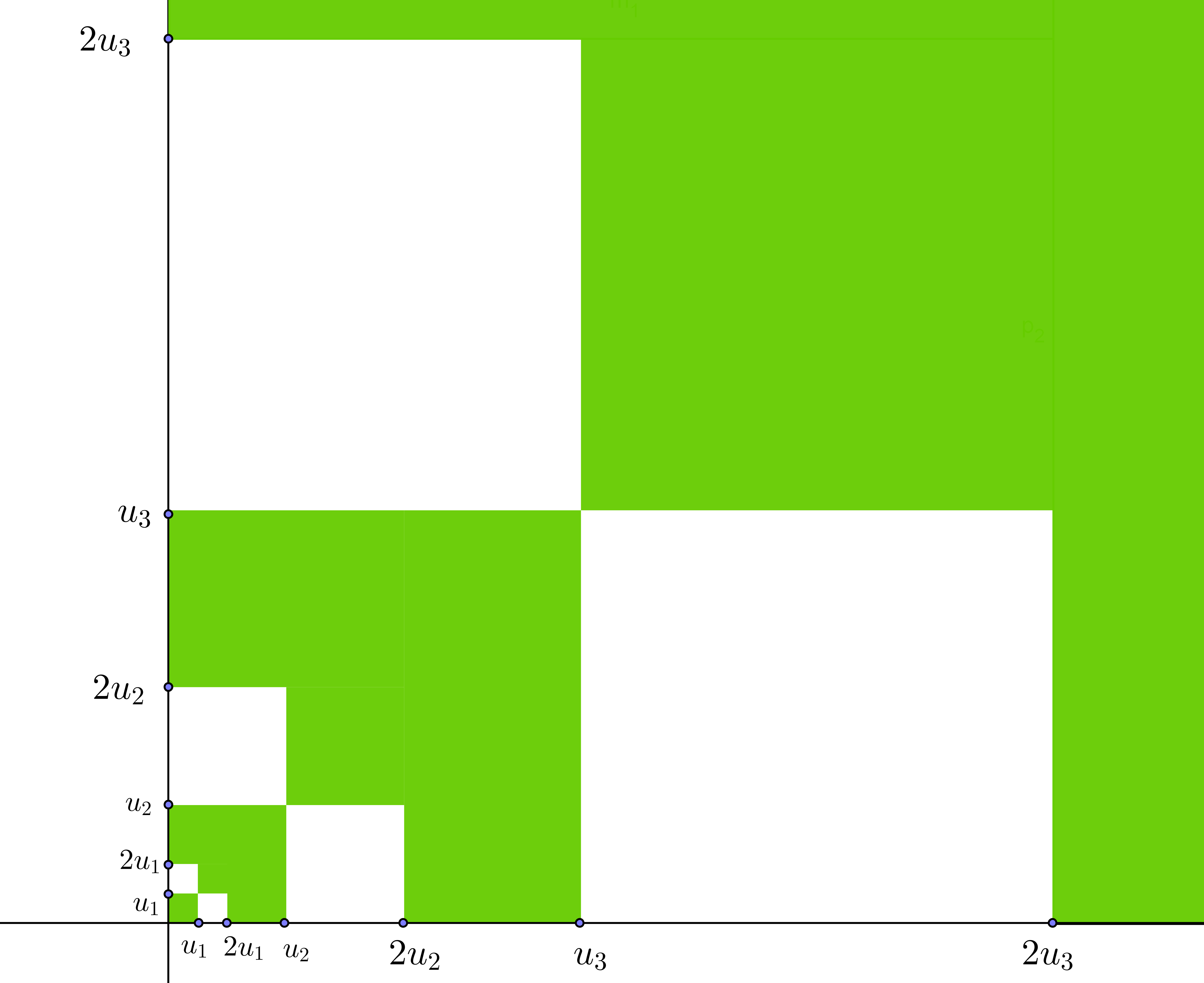}
\centering
\caption{The set $A$ in Example~\ref{eg: fractal set example}.}
\label{Fig: fractal example}

\end{figure} One sees that $\underline{d}_{\textbf{Rect}}(A)= \frac{1}{2}$ and $\overline{d}_{\textbf{Rect}}(A)=1$ (for details, see Proposition\ref{formula for density of fractal}). However one may still give a lower bound for $\underline{d}_{\textbf{Rect}}(A+B)$ (which is greater than the trivial lower bound $\underline{d}_{\textbf{Rec}}(A)$) for arbitrary $B \subset \mathbb{N}^2$ such that $B+B=\mathbb{N}^2$, as follows. It is not hard to see (Proposition~\ref{formula for density of fractal}) that $\underline{d}_{\textbf{Tab}}(A) =\frac{1}{3}$. Thus by Theorem~\ref{thm: plunnecke for d_tab} we have that \begin{align*} \underline{d}_{\textbf{Rect}}(A+B) \geq \underline{d}_{\textbf{Tab}}(A+B) \geq \underline{d}_{\textbf{Tab}}(A)^{1/2} \underline{d}_{\textbf{Tab}}(B+B)^{1/2}=\frac{1}{\sqrt{3}}. \end{align*}  

\end{eg}

One can generalise the example above to construct more general sets that are asymptotically finite unions of translates of a large square. 

\begin{mydef}[Fractal set generated by a pattern] \label{def: fractal set} As before, we will use the convention $[a,b) = \{x \in \mathbb{N}\text{ } |\text{ } a \leq x <b \}$. Let $N \in \mathbb{N}$ and $$P \subset \{0,1, \ldots, N\} \times \{0,1,\ldots, N \} $$ be a set such that $(0,0) \in P$. We call such a $P$ a \textit{pattern}. Choose a sequence $u_1,u_2, \ldots$ of positive integers such that $u_k>(N+1)u_{k-1}$ and $\lim_{k \to \infty} \frac{u_{k}}{u_{k-1}} = \infty$. Define $$\cP_k = \{ u_k p \text{ }| \text{ } p \in P \} + [0,u_k)^2$$ and let $$A_k = \cP_k \setminus [0,(N+1)u_{k-1})^2$$ where $u_0=0$ (in other words $A_1=\cP_1$). Finally, we define $A(P,N) = \bigcup_{k=1}^{\infty} A_k$ to be \textit{the fractal set generated by the pattern $P$ of degree $N$}.

\end{mydef}

\begin{figure}[H]
\includegraphics[scale=0.6]{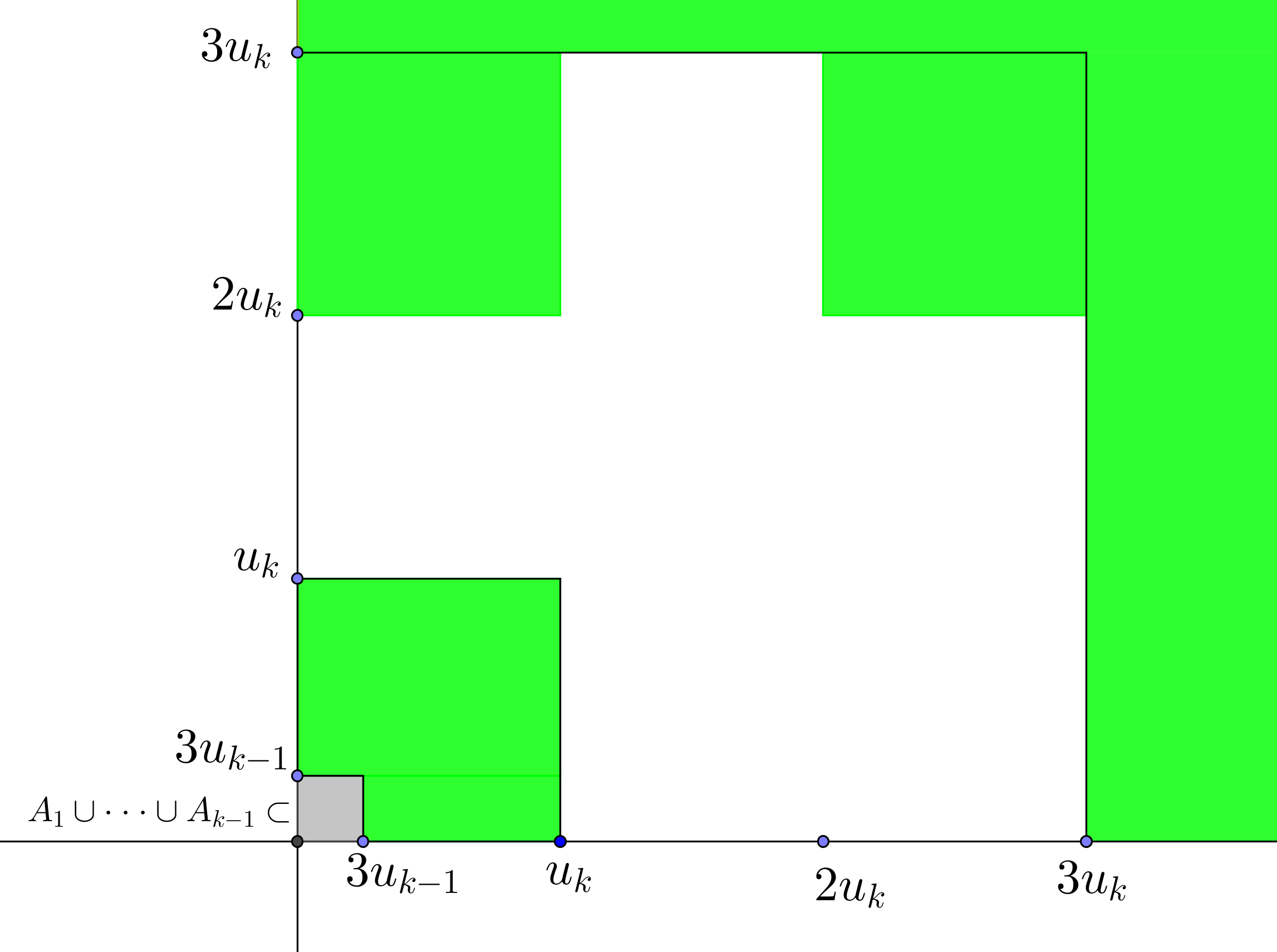}
\centering
\caption{Definition~\ref{def: fractal set} for the case $N=2$, $P=\{(0,0),(0,2),(2,2)\}$. Note that $A_k$ is contained in $[0,3u_k)^2 \setminus [0,3u_{k-1})^2$; while $\cP_{j}$ and $A_{j}$ are contained in the small square $[0,3u_{k-1})^2$, for all $j<k$.}
\label{Fig: fractal recursion example}

\end{figure}

Note that the set $A$ given in Example~\ref{eg: fractal set example} is $A(P,1)$ for $P = \{ (0,0),(1,1) \}$. We now give a combinatorial formula for $\underline{d}_{\textbf{Tab}}$ and $\underline{d}_{\textbf{Rect}}$ of fractal sets generated by a pattern.

\begin{prop}\label{formula for density of fractal} Suppose that $N \in \mathbb{N}$ and $P \subset \{0,1, \ldots, N\} \times \{0,1,\ldots, N\}$ and let $A=A(P,N)$. Then $$\underline{d}_{\textbf{Rect}}(A)=\min \left\{ \frac{|P \cap I|}{|I|} \text{ } | \text{ } I=[0,m] \times [0,n] \cap \mathbb{N}^2 \text{ for some } m,n \in \{0,1, \ldots, N\}  \right \}$$ and $$\underline{d}_{\textbf{Tab}}{A} = \min \left \{ \frac{|P \cap T|}{|T|} \text{ } | \text{ } T \subset \{0, \ldots, N\} \times \{0, \ldots, M\} \text{ is a tableau.}  \right \}$$ where we define a \textit{tableau} to be a set of the form $$\bigcup_{(n,m) \in F} \{0,1, \ldots, n\} \times \{0,1,\ldots,m\}$$ for some finite $F \subset \mathbb{N}^2$.

\end{prop}

\textbf{Proof:} See Appendix~\ref{section: Densities of fractal sets}. $\blacksquare$

Thus  $\underline{d}_{\textbf{Rect}}(A) < 1= \overline{d}_{\textbf{Rect}}(A)$ for all fractal sets other than $\mathbb{N}^2$, which leads to many applications of our results that do not follow from Proposition~\ref{prop: plunnecke when d exists}.

\begin{eg} Let $P=\{(0,0),(0,2),(2,2)\}$ and let $A=A(P,2)$ (as depicted in Figure~\ref{Fig: fractal recursion example}). Then $$\underline{d}_{\textbf{Rect}}(A)=\underline{d}_{\textbf{Tab}}(A) = \frac{1}{6}.$$ Thus if $(0,0) \in B \subset \mathbb{N}^2$ and $k>0$ is an integer then by Theorem~\ref{thm: drect bound in terms of dtab} we have $$\underline{d}_{\textbf{Rect}}(A+k'B) \geq \underline{d}_{\textbf{Rect}}(A)^{1-\frac{k'}{k}} \underline{d}_{\textbf{Rect}}(kB)^{\frac{k'}{k}}$$ for integers $0<k'<k$. Many such special cases of Question~\ref{question: plunnecke for d_rect} may be constructed, which are not covered by Proposition~\ref{prop: plunnecke when d exists}.

\end{eg}

\subsection{A remark about higher dimensions} When $B$ is an asymptotic basis of order $k$ (i.e., $kB=\bN^2$), the natural analogue of Theorem~\ref{main result in intro} (and hence also its consequences) holds in higher dimensions, as the reader may observe in the equation (\ref{special case when kB=N^2}) found in the proof of Theorem~\ref{main result in intro} in Section~\ref{section: proof of main plunnecke inequality}. However, it is yet unclear whether the higher dimensional analogue of Theorem~\ref{main result in intro} holds in full generality, the main obstacle lies in finding the right higher dimensional extension of Lemma~\ref{Lemma: constructing a_j}.

\subsection{Organization of the paper}

In Section~\ref{sec: Plunnecke inequalities for truncated sumsets} we review some classical Pl\"unnecke inequalities for cardinalities of truncated sumsets (i.e., sets of the form $(A+B) \setminus C$), as well as the less well known but crucial $\delta$-heavy Pl\"unnecke inequality. In Section~\ref{sec: Tableaux} we introduce the main definitions, notations and conventions. In particular, we introduce the notion of a \textit{tableau}, the main combinatorial object in this paper, and we establish some basic properties of the lower tableau density $\underline{d}_{\textbf{Tab}}$. Sections~\ref{section: Trimming lemma} and \ref{section: approximating subtableaux} include some technical combinatorial lemmata involving tableaux that will be put together in Section~\ref{section: proof of main plunnecke inequality} to conclude the proof of Theorem~\ref{main result in intro}. In Section~\ref{sec: questions} we state some related open problems; more specifically, some conjectural multidimensional Pl\"unnecke inequalities for various densities. Finally, Appendix~\ref{section: Densities of fractal sets} is devoted to proving Proposition~\ref{formula for density of fractal}.

\textbf{Acknowledgement:} The author is grateful for many insightful and encouraging conversations with Alexander Fish.

\section{Pl\"unnecke inequalities for truncated sumsets}

\label{sec: Plunnecke inequalities for truncated sumsets}

One of the key tools used in the proof of our results, as well as Jin's proofs in \cite{JinPlunnecke} and \cite{JinEpsilon}, is the following Pl\"unnecke inequality for truncated sumsets.

\begin{thm}[See \cite{Ruzsasumsetsandstructure}] \label{thm: truncated sumsets} Let $A,B,C$ be finite subsets of an abelian group and define $$D_n=D_n(A,B,C)=\min_{\emptyset \neq A' \subset A} \frac{|(A'+nB) \setminus (C+(n-1)B)|}{|A'|}. $$ Then $D^{1/n}_n$ is decreasing in $n$.

\end{thm}

We will need (unlike Jin in \cite{JinPlunnecke} and \cite{JinEpsilon}) the following $\delta$-heavy version of this inequality. We will include a proof for the sake of completeness as this version, to the best of the author's knowledge, rarely appears in the literature (cf. \cite{BjorklundFishPlunnecke}).

\begin{thm}\label{A B C inequality}  Let $A,B,C$ be finite subsets of an abelian group and let $0<\delta<1$. Then for positive integers $k' <k$, there exists $A' \subset A$ such that $|A'| > \delta|A|$ and $$\frac{|(A'+kB) \setminus (C+(k-1)B) |}{|A'|}  \leq (1-\delta)^{-k/k'} \left( \frac{|(A+k'B) \setminus (C+(k'-1)B)|}{|A|} \right)^{k/k'}.  $$ 

\end{thm}

\textbf{Proof:} Using Theorem~\ref{thm: truncated sumsets} and the fact that $(1-\delta)^{-k/k'}>1$, take non-empty $A' \subset A$ of maximal cardinality such that \begin{align}\label{eq1 in delta proof} |(A'+kB) \setminus (C+(k-1)B) |  \leq (1-\delta)^{-k/k'} \left( \frac{|(A+k'B) \setminus (C+(k'-1)B)|}{|A|} \right)^{k/k'} |A'|. \end{align} Suppose for contradiction that $|A'| \leq \delta|A|$, thus $|A \setminus A'| \geq (1-\delta)|A|$. Now apply Theorem~\ref{thm: truncated sumsets}, with $A \setminus A'$ playing the role of $A$, to obtain a non-empty $A'' \subset A \setminus A'$ such that 

$$|(A''+kB) \setminus (C+(k-1)B) |  \leq \left( \frac{|((A\setminus A')+k'B) \setminus (C+(k'-1)B)|}{|A\setminus A'|} \right)^{k/k'} |A''|.$$

We deduce from this inequality, using the fact that $A \setminus A' \subset A$ and $|A \setminus A'| \geq (1-\delta)|A|$, the estimate \begin{align}\label{eq2 in delta proof} |(A''+kB) \setminus (C+(k-1)B) |  \leq (1-\delta)^{-k/k'} \left( \frac{|(A+k'B) \setminus (C+(k'-1)B)|}{|A|} \right)^{k/k'} |A''| . \end{align} Now by adding the estimates (\ref{eq1 in delta proof}) and (\ref{eq2 in delta proof}) we get $$|((A' \cup A'') + kB) \setminus (C+(k-1)B)| \leq (1-\delta)^{-k/k'} \left( \frac{|(A+k'B) \setminus (C+(k'-1)B)|}{|A|} \right)^{k/k'} |A' \cup A''|.  $$ But since $A''$ is non-empty, this contradicts the maximality of $A'$. $\blacksquare$ 

\section{Tableaux}

\label{sec: Tableaux}

We stress that throughout this paper we use (and already have used) the convention $\mathbb{N}=\{0,1, \ldots \}$. A \textit{tableau} is a set of the form $$ \bigcup_{(N,M) \in F} \{0, 1, \ldots, N\} \times \{0,1,\ldots M\} $$ where $F \subset \mathbb{N}^2$ is a finite set. It will also be convenient to define a \textbf{tableau region} to be a set of the form $$ \bigcup_{(N,M) \in F} [0,N) \times [0,M) $$ for some finite $F \subset \mathbb{N}^2$. Thus a tableau is precisely the set of lattice points of some tableau region.

\textbf{Important note on notation:} If $A$ is a union of rectangles of the form $[N,N+1) \times [M,M+1)$ where $N,M \in \bZ$ (such as a tableau region), then by $|A|$ we mean the Lebesgue measure of $A$, which is also the number of integer points in $A$. For most of our arguments, it is more conceptual to consider the more geometric and continuous notion of Lebesgue measure. Unless otherwise specified, by $[a,b)$ we mean $\{x \in \bR \text{ }|\text{ } a \leq x < b \}$ (very rarely it will mean $\{x \in \bZ \text{ }|\text{ } a \leq x < b \}$, in fact only in the definition of \textit{fractal sets} in the Introduction above and Appendix~\ref{section: Densities of fractal sets}).

We note the following simple but useful additive characterization of tableaux.

\begin{lemma} If $T \subset \mathbb{N}^2$ is finite and non-empty, then the following are equivalent: 
\begin{enumerate}[label=({\roman*})]
	\item $T$ is a tableau.
	\item $\mathbb{N}^2 \setminus T$ is invariant under addition by elements of $\mathbb{N}^2$, i.e., $(\bN^2 \setminus T)+a \subset \bN^2 \setminus T$ for all $a \in \bN^2$.
\end{enumerate} 
\end{lemma}

This means that if $T$ is a tableau and $B \subset \bN^2$ contains $(0,0)$, then $(\bN^2 \setminus T) + B =\bN^2 \setminus T$. As a consequence, we may apply Theorem~\ref{A B C inequality} with $C=\bN^2 \setminus T$ to obtain the following crucial proposition. 

\begin{prop}\label{truncated heavy plunnecke} Let $T \subset \mathbb{N}$ be a tableau and suppose that $A,B \subset \mathbb{N}^2$ with $(0,0) \in B$. Then for positive integers $k' <k$ and $0<\delta<1$, there exists $A' \subset A$ such that $|A'| > \delta|A|$ and $$\frac{|(A'+kB) \cap T|}{|A'|}  \leq (1-\delta)^{-k/k'} \left( \frac{|(A+k'B) \cap T|}{|A|} \right)^{k/k'}.  $$ 

\end{prop}

We now turn to explicating some basic properties of the densities $\underline{d}_{\textbf{Tab}}$ and $\underline{d}_{\textbf{Tab}(L)}$ that were introduced above. The following simple lemma will be convenient as it shows that we may, without loss of generality, assume that the side lengths of our rectangles are divisible by a chosen integer.

\begin{lemma} \label{Lemma: sidelengths divide D}Let $D,L$ be positive integers and $A \subset \mathbb{N}^2$. Then there exists a sequence of the form $$ F_n=\bigcup_{j=1}^{\ell(n)} [0, W_{j,n}) \times [0, H_{j,n})$$ with each $\ell(n) \leq L$, such that

\begin{enumerate}[label=(\textbf{\roman*})] 
	\item For each $j \in \{1,\ldots, \ell(n) \}$, we have $$W_{j,n} \equiv H_{j,n} \equiv 0 \mod D.$$
	\item The $W_{i,n}$ and $H_{i,n}$ tend to $\infty$, more precisely $$\lim_{n \to \infty} \min_{i \in \{1, \ldots, \ell(n)\}} W_{i,n} = \lim_{n \to \infty} \min_{i \in \{1, \ldots, \ell(n)\}} H_{i,n}= \infty.$$
	\item \label{condition: approx density} \begin{align*} \liminf_{n \to \infty} \frac{|A \cap F_n|}{|F_n|} = \underline{d}_{\textbf{Tab}(L)}(A). \end{align*} 

\end{enumerate}

\textbf{Proof Sketch:} The idea is that if one replaces $W_{i,n}$ with $W_{i,n} + O(1)$ then $$\liminf_{n \to \infty} \frac{|A \cap F_n|}{|F_n|} $$ remains unchanged. Applying this finitely many times, we can adjust a sequence satisfying \ref{condition: approx density} (which exists by a simple diagonalization argument) to one which satisfies the desired properties. $\blacksquare$

\end{lemma}

Let us spell out a useful characterization of the lower tableaux density. 

\begin{lemma}\label{characterization of lower tableaux density} Let $A \subset \mathbb{N}^2$ with $\alpha=\underline{d}_{\textbf{Tab}}(A)$. Then for each $\epsilon>0$ and positive integer $L$, there exists an $R=R_{\epsilon,L}>0$ such that whenever $M_{i}, N_{i}>R$ are integers, for $i = 1,2, \ldots L$, we have $$\frac{|A \cap F|}{|F|} > \alpha -\epsilon$$ where $$F=\bigcup_{i=1}^{L} [0,N_{i}) \times [0,M_{i}) \cap \mathbb{Z}^2.$$ Moreover, $\underline{d}_{\textbf{Tab}}(A)$ is the largest choice of $\alpha$ which makes this statement true.

\end{lemma}

\section{Trimming lemma}

\label{section: Trimming lemma}

In this section we formulate and prove the \textit{Trimming Lemma}, one of the main combinatorial tricks of this paper. It will be most convenient to state and prove it in a rather abstract setting. If $X$ is a set equipped with a measure $\mu$, then we will use the \textit{averaging notation} $$A_{\mu}(f,U)=\frac{1}{\mu(U)}\int_{U} f d\mu$$ for $f:X \to \mathbb{R}$ and $U \subset X$. If the measure $\mu$ is clear, we simply use the shorthand $A(f,X)$.

\begin{lemma}\label{trimming lemma} Let $I \subset \mathbb{N}^2$ be a tableau equipped with a positive measure $\mu$ (on the set of all subsets of $I$). Suppose that $$\rho: I \to [0,1]$$ is a function and $\alpha>0$ is such that $$ A(\rho,S) \geq \alpha$$ for all non-empty tableaux $S \subset I$. Then there exists $$\rho': I \to [0,1]$$ such that 

\begin{enumerate}[label=(\textbf{\roman*})]
	\item $\rho' \leq \rho$
	\item \label{overall not too skinny}$A(\rho',I) \geq \alpha$
	\item \label{upper regions trimmed} For all tableaux $S \subsetneqq I$, we have $$A(\rho',I \setminus S) \leq \alpha.$$
\end{enumerate}

\end{lemma}

\begin{figure}[H]
\centering
\includegraphics[scale=1.0]{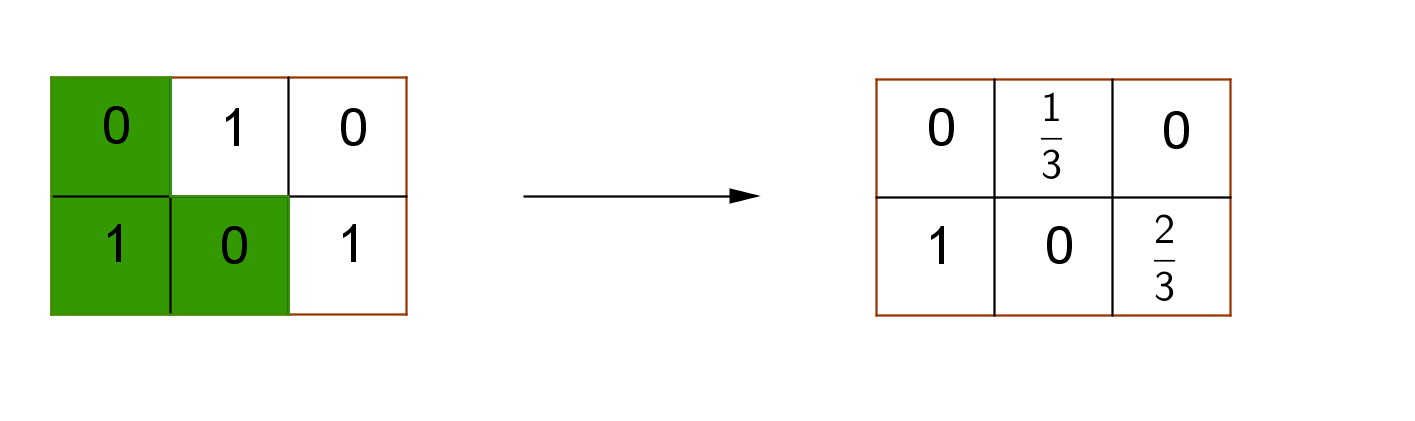}
\caption{Example of the theorem. The shaded tableaux shows that $\alpha=1/3$ is the maximal choice of $\alpha$ in the hypothesis. The rightmost column has average value of $\rho$ equal to $1/2$, so we must have $\rho' \neq \rho$ for all suitable $\rho'$.}
\label{Fig: example of thm}

\end{figure}

\begin{eg} Let $I=\{0,1,2\} \times \{0,1\}$, $\mu$ be the counting measure and $\rho$ be given as on the left of Figure~\ref{Fig: example of thm}. A choice of $\alpha = \frac{1}{3}$ satisfies the hypothesis of the theorem. In fact, the shaded tableaux shows that it is the largest choice of $\alpha$. On the right we have a possible choice of $\rho'$.

\end{eg}

\begin{remark} Equality actually occurs in \ref{overall not too skinny} since $S$ may taken to be empty in \ref{upper regions trimmed}.\end{remark}

\textbf{Proof of Lemma~\ref{trimming lemma}:} We proceed by induction on $|I|$. The case $I=\{(0,0)\}$ is clear, one can just set $\rho'((0,0))=\alpha$. Now suppose the theorem holds for all tableaux with cardinality strictly less than $|I|$. Suppose that \ref{upper regions trimmed} fails for some tableau $S \subsetneqq I$ with $\rho'=\rho$ (otherwise, we may take $\rho'=\rho$). Let $S_{max} \subsetneqq I$ be a maximal tableau contained in $I$ such that $$A(\rho,I \setminus S_{max}) > \alpha.$$ For $x \in I \setminus S_{max}$ we define \begin{align} \label{definition of rho' outside of Smax} \rho'(x)=\frac{\alpha}{A(\rho,I \setminus S_{max})} \rho(x). \end{align}

To define $\rho'$ on $S_{max}$, we use the induction hypothesis as follows. Let $\rho_0=\rho|_{S_{max}}$. Since $|S_{max}|<|I|$ we apply the induction hypothesis to $\rho_0$ to obtain a map $\rho_0':S_{max} \to [0,1]$ such that 

\begin{enumerate}[label=(\textbf{\alph*})]
	\item $\rho'_0 \leq \rho_0$.
	\item \label{induction: overall not too skinny} $A(\rho'_0, S_{max}) \geq \alpha.$
	\item \label{induction: upper region trimmed} For all tableau $S \subsetneqq S_{max}$ we have $$A(\rho'_0, S_{max} \setminus S) \leq \alpha. $$

\end{enumerate}

We define $\rho'(x)=\rho'_0(x)$ for $x \in S_{max}$. Let us now check that $\rho':I \to [0,1]$ satisfies the desired conclusions. We have $\rho' \leq \rho$ since $\rho_0' \leq \rho_0$ and $$\frac{\alpha}{A(\rho, I \setminus S_{max})} < 1.$$ By (\ref{definition of rho' outside of Smax}) we have $$A(\rho', I \setminus S_{max}) = \alpha$$ which together with \ref{induction: overall not too skinny} implies that \ref{overall not too skinny} holds. Now suppose that $S \subsetneqq I$ is a tableau. Then we may decompose $$I \setminus S= (S_{max} \setminus (S \cap S_{max})) \bigsqcup (I \setminus (S_{\max} \cup S)). $$ It is enough to show that $$A(\rho', X) \leq \alpha$$ when $X$ is one of these parts. If $$S_{\max} \setminus (S \cap S_{\max}) \neq \emptyset$$ then, by \ref{induction: upper region trimmed}, we do indeed have $$A(\rho', S_{\max} \setminus (S \cap S_{\max})) \leq \alpha$$ since $S \cap S_{\max} \subset S_{\max}$ is a tableau. Now suppose that $$I \setminus (S_{\max} \cup S) \neq \emptyset$$ and consider the following two cases.

\underline{\textbf{Case 1:} $I \setminus (S_{\max} \cup S) = I \setminus S_{\max}$ }

Then we have $$A(\rho', I \setminus (S_{\max} \cup S))= A(\rho', I \setminus S_{\max}) = \frac{\alpha}{A(\rho, I \setminus S_{\max})} A(\rho, I \setminus S_{\max}) = \alpha. $$

\underline{\textbf{Case 2:} $I \setminus (S_{\max} \cup S) \subsetneqq I \setminus S_{\max}$}

This means that $S_{\max} \subsetneqq S_{max} \cup S$ and thus, by the maximality of $S_{\max}$, we have that $$A(\rho', I \setminus (S_{max} \cup S)) \leq A(\rho, I \setminus (S_{max} \cup S)) \leq \alpha.$$

This verifies \ref{upper regions trimmed} and thus completes the proof. $\blacksquare$

\section{$Q^2$-tilings and approximating subtableaux regions}

We now turn to studying tableau regions obtained by subdividing a tableau region. We will consider subdivisions that are equally spaced, thus it is convenient to define following the notion.

\label{section: approximating subtableaux}

\begin{mydef}[$D$-tableau region] If $D$ is a positive integer, then we define a \textit{$D$-tableau} region to be a tableau region where all side lengths are divisible by $D$. More precisely, a $D$-tableau is a set of the form $$\bigcup_{(N,M) \in F} [0,N) \times [0, M) $$ for some finite $F \subset D\bN^2$.

\end{mydef}

\begin{lemma}\label{lemma: approximating grid} Fix a positive integer $Q$. Let $F=[0,N] \times [0,M]$, where $N,M$ are positive integers divisible by $Q$, and suppose that $S=F \setminus S'$ for some tableaux region $S' \subset F$. Let $$\mathcal{C}_0=\left \{ \left[\frac{Ni}{Q},\frac{N(i+1)}{Q} \right) \times \left[ \frac{Mj}{Q}, \frac{M(j+1)}{Q} \right) | i,j \in \{0,1\ldots, Q-1\} \right\}.$$ Let $$\widehat{S} = \bigcup_{C \in \mathcal{C}_0, C \cap S \neq \emptyset} C$$be the smallest set that contains $S$ and is in the $\sigma$-algebra generated by the partition $\mathcal{C}_0$. Then $$\frac{|\widehat{S}| - |S|}{NM} \leq \frac{2}{Q}.$$

\end{lemma}

\begin{figure}[H]
\centering
\includegraphics[scale=0.20]{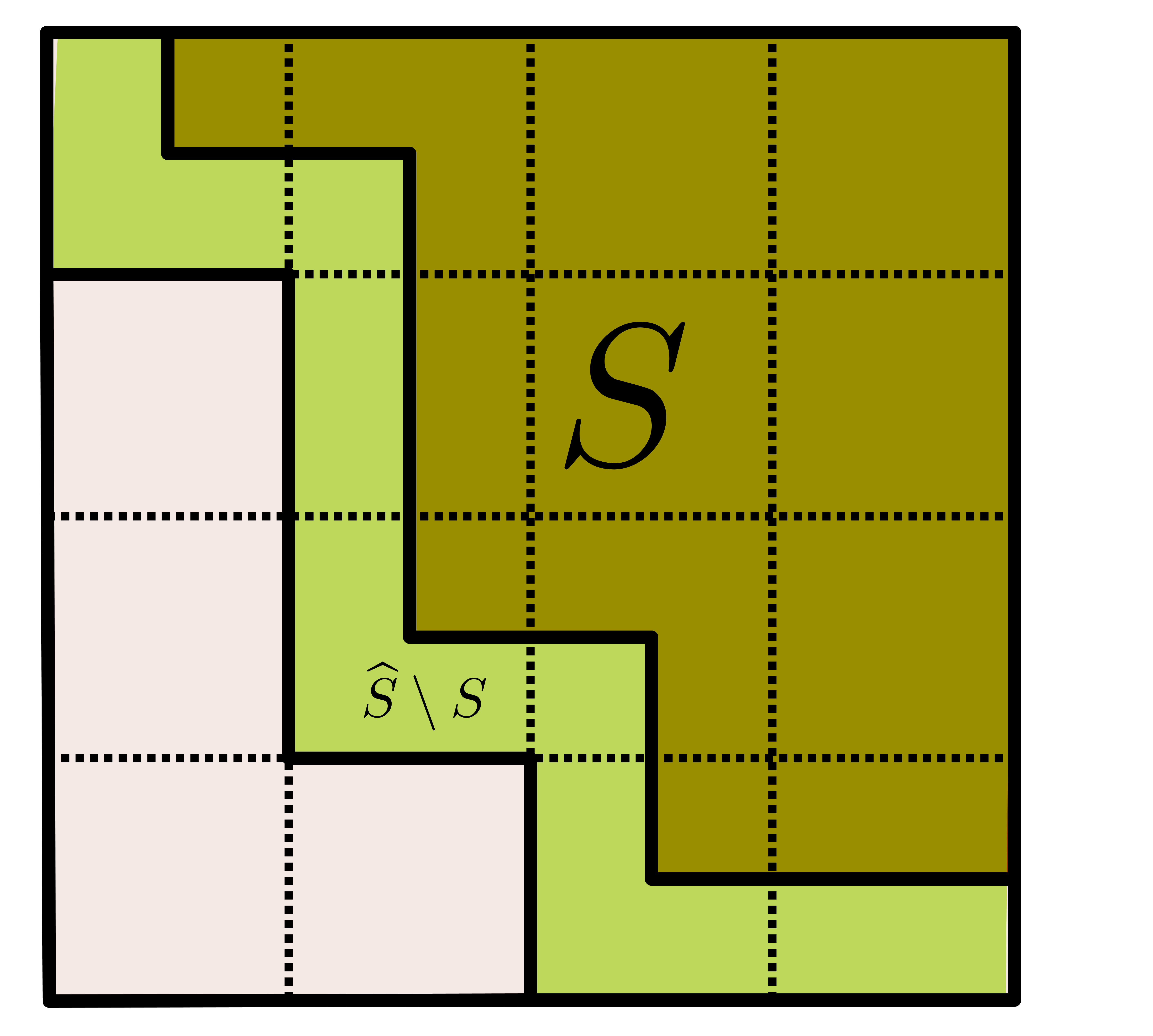}
\caption{ An example, with $Q=4$ of Lemma~\ref{lemma: approximating grid}. We can cover $\widehat{S} \setminus S$ with $7$ of the $16$ square tiles. }
\label{fig: basic_rectangle_estimate}
\end{figure}

\textbf{Proof:} There are at most $2Q$ elements of $\mathcal{C}$ that intersect both $S$ and $F \setminus S$, since these elements form a path consisting of right and down steps. In fact, there are at most $2Q-1$. See Figure~\ref{fig: basic_rectangle_estimate}. $\blacksquare$

\begin{mydef}($Q^2$-tiling and its refinement) \label{def: Q^2 tiling} Let $Q$ be a positive integer and suppose that $$F= \bigcup_{m=1}^{\ell} [0,W_m) \times [0,H_m)$$ with $$H_1>H_2 > \ldots > H_{\ell}>0,$$ $$0<W_1<W_2<\ldots <W_{\ell}$$ and $$W_m \equiv H_m \equiv 0 \mod Q^2$$ for all $m \in \{1, \ldots, \ell\}$. Then we define the \textbf{$Q^2$-tiling} $$\mathcal{C}_0=\mathcal{C}_0(F,Q) = \left\{ C^m_{i,j} \text{ } | \text{ } i,j \in \{0, \ldots, Q-1\} \text{ and } m \in \{1, \ldots, \ell\} \right \} $$ where $$C^m_{i,j} = \left [ W_{m-1} + \frac{i}{Q}(W_{m}-W_{m-1}), W_{m-1} + \frac{i+1}{Q}(W_{m}-W_{m-1}) \right) \times \left[\frac{j}{Q}H_m, \frac{j+1}{Q}H_m \right) $$ where $W_0=0$. We may refine $\mathcal{C}_0$ as follows: Let $\mathcal{Y}=\mathcal{Y}(F,Q)$ be the set of integers which appear as a $y$ ordinate of some corner of a cell $C^m_{i,j}$, in other words $$\mathcal{Y}= \left \{ \frac{j}{Q}H_m \text{ }|\text{ } m \in \{1, \ldots, \ell\}, j \in \{0,1, \ldots, Q\} \right\}.$$ Now write $$\mathcal{Y}=\{0=y_1 < y_2 < \ldots < y_{|\mathcal{Y}|}\}$$ and let $\mathcal{H}= \{F \cap \left( \bR \times [y_r,y_{r+1}) \right) \text{ } | \text{ } r \in \{0,1, \ldots, Q-1 \} \}$ which is a partition of $F$. Define the \textbf{refined $Q^2$-tiling} of $F$ to be the common refinement $$\mathcal{C}=\mathcal{C}(F,Q)=\mathcal{C}_0 \lor \mathcal{H} $$ of the partitions $\mathcal{C}_0$ and $\cH$ of $F$. More explicitly, $$\mathcal{C}=\mathcal{C}(F,Q)=\left\{ C^m_{i,j} \cap \mathbb{R} \times [y_{r},y_{r+1}) \text{ }| \text{ } C^m_{i,j} \in \mathcal{C}_0 \text{ and } r \in \{1, \ldots, |\mathcal{Y}|-1\} \right\} \setminus \{ \emptyset \}.$$ 

\begin{figure}[H]
\centering
\includegraphics[scale=0.21]{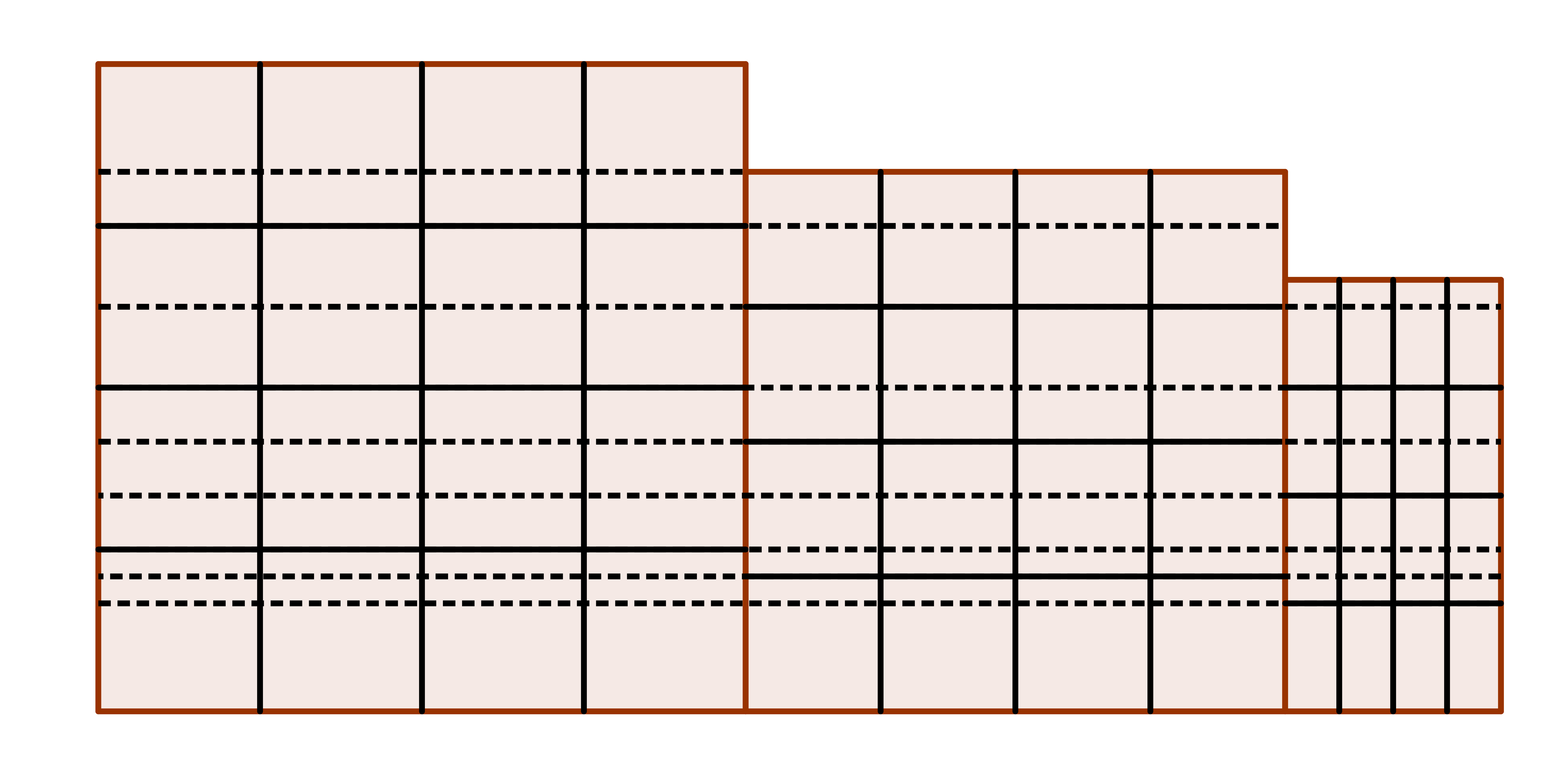}
\caption{ An example, with $Q=4$, of refining $\mathcal{C}_0$ to $\mathcal{C}$.}
\label{fig: C0 tiling example}
\end{figure}

Note that the elements of $\mathcal{C}$ are rectangles with both side lengths integers divisible by $Q$, and thus contain at least $Q^2$ elements of $\mathbb{N}^2$. 

Naturally, we may identify $\mathcal{C}$ with a tableau $T=T(F,Q) \subset \mathbb{N}^2$ by constructing a bijection $\psi:\mathcal{C} \to T$ as follows: \begin{enumerate}[label=(\roman*)]
	\item $\psi^{-1}(0,0)$ is the unique element of $\mathcal{C}$ which contains $(0,0)$.
	\item  $\psi^{-1}(i+1,j)$ is the element to the right of $\psi^{-1}(i,j)$ and $\psi^{-1}(i,j+1)$ is the element just above $\psi^{-1}(i,j)$. (Note: We refined $\mathcal{C}_0$ to $\mathcal{C}$ precisely so that the notion of \textit{right} is well defined.)
\end{enumerate}

A set of the form $\bigcup \psi^{-1}(T') $ (we use the notation $\bigcup X=\bigcup_{x \in X} x$), for some tableau $T' \subset T$, will be called a \textbf{$\mathcal{C}$-measurable subtableau region.} In general, a union of elements of $\cC$ will be called a $\cC$-measurable set.

\begin{remark}\label{remark about sidelengths} All $\mathcal{C}$-measurable subtableau regions contain the element of $\mathcal{C}$ that contains $(0,0)$. This element is precisely $$[0, \frac{1}{Q}W_1) \times [0, \frac{1}{Q}H_{\ell}).$$ Hence, each $\mathcal{C}$-measurable subtableau region is a union of rectangles of width at least $\frac{1}{Q}W_1$ and height $\frac{1}{Q}H_{\ell}$.

\end{remark}

\end{mydef}

\begin{lemma}($\mathcal{C}$-measurable approximations)\label{lemma: approximating regions} Fix a positive integer $Q$ and a $Q^2$-tableau region $F$. Suppose that $S \subset F$ is of the form $$S=F \setminus F'$$ where $F' \subset F$ is a tableau region. Let $$\widetilde{S} = \bigcup_{C \in \mathcal{C}, C \cap S \neq \emptyset} C$$ be the smallest subset of $F$ that contains $S$ and may be written as a union of elements of $\mathcal{C}=\mathcal{C}(F,Q)$. Then $$\frac{|\widetilde{S} \setminus S|}{|F|} \leq \frac{2}{Q}. $$

\end{lemma} 

\textbf{Proof:} We will use the setup from Definition~\ref{def: Q^2 tiling} above (i.e. the parameters $m$, $H_i$, $W_i$ etc.). For each $m \in \{1, \ldots, \ell \}$ define $$U_m=[W_{m-1},W_m) \times [0,H_m) = \bigcup_{(i,j) \in \{0,1, \ldots Q-1 \}^2} C^m_{i,j}. $$ Apply Lemma~\ref{lemma: approximating grid} to $S \cap U_m \subset U_m$ to get $$|(\widetilde{S} \setminus S) \cap U_m | \leq |(\widehat{S} \setminus S) \cap U_m| \leq \frac{2|U_m|}{Q}$$ where $\widehat{S} \supset \widetilde{S}$ is the smallest subset of $F$ that contains $S$ and may be written as a union of elements of $\mathcal{C}_0$. Since $F=\bigsqcup_{m=1}^{\ell} U_m$ we are done by summing this estimate over $m \in \{1, \ldots, \ell\}$. $\blacksquare$

\begin{lemma}[Trimming a set of points]\label{lemma: trimming set of points} Fix a positive integer $Q$, a $Q^2$-tableaux region $F$ and $A \subset F \cap \mathbb{N}^2$. Let $\mathcal{C}=\mathcal{C}(F,Q)$. Define $$\alpha = \inf \left\{ \frac{|A \cap F'|}{|F'|}\text{ }| \text{ }F' \subset F \text{ is a non-empty $\mathcal{C}$-measurable subtableau region} \right\}.$$ Then there exists $A' \subset A$ such that 

\begin{enumerate}[label=(\alph*)]
	\item For $\mathcal{C}$-measurable subtableau regions $F' \subsetneqq F$ we have $$ \frac{|A' \cap (F \setminus F')|}{|F \setminus F'|} \leq \alpha + \frac{1}{Q^2}.$$
	\item $$\frac{|A' \cap F|}{|F|} \geq \alpha.$$
\end{enumerate}

\end{lemma}

\textbf{Proof:} $T=T(F,Q) \subset \mathbb{N}^2$ be the corresponding subtableau and let $\psi:\mathcal{C} \to T$ be the bijection constructed in Definition~\ref{def: Q^2 tiling}. Apply the Trimming Lemma (Lemma~\ref{trimming lemma}) to the tableaux $T$, the measure given by $\mu(\{t\}) = |\psi^{-1}(t)|$ for $t \in T$ and the map $\rho:T \to [0,1]$ given by $$ \rho(t)= \frac{|A \cap \psi^{-1}(t)|}{|\psi^{-1}(t)|}$$ to obtain $\rho':T \to [0,1]$ such that 

\begin{enumerate}[label=(\roman*)]
	\item $\rho' \leq \rho$
	\item For all tableaux $T' \subsetneqq T$, we have $$\frac{1}{|\bigcup \psi^{-1}(T \setminus T')|} \sum_{t \in T \setminus T'} |\psi^{-1}(t)|\rho'(t) \leq \alpha.$$ 
	\item $$\frac{1}{|F|}\sum_{t \in T}|\psi^{-1}(t)|\rho'(t) \geq \alpha.$$
\end{enumerate}

Since each element of $\mathcal{C}$ is a rectangle with both sidelengths multiples of $Q$, we may find $A' \subset A$ such that, for all $t \in T$, $$\rho'(t) \leq \frac{|A' \cap \psi^{-1}(t)|}{|\psi^{-1}(t)|} \leq \rho'(t) + \frac{1}{Q^2}.  $$ $\blacksquare$

We use the notation $y(a,b)=b$. 

\begin{lemma} \label{Lemma: constructing a_j} Fix a positive integer $Q$, a $Q^2$-tableaux region $F$ and $A \subset F \cap \mathbb{N}^2$ and let $$S= \bigcup_{a \in A} (a + [0,\infty)^2) \cap F.$$ Then there exists a positive integer $J$ and a sequence $a_1,a_2, \ldots, a_J \in A$ such that $$y(a_j) \leq y(a_{j-1})-Q$$ and $$G := \bigcup_{j=1}^J (a_j + [0,\infty)^2) \cap F$$ satisfies $$\frac{|S| - |G|}{|F|} \leq \frac{3}{Q}.$$

\end{lemma}

\textbf{Proof:} Note that $S=F \setminus F'$ for some tableau region $F'$. As in Lemma~\ref{lemma: approximating regions} we let $$\widetilde{S} = \bigcup_{C \in \mathcal{C}, C \cap S \neq \emptyset} C$$ be the smallest subset of $F$ that contains $S$ and may be written as a union of elements of $\mathcal{C}=\mathcal{C}(F,Q)$. So $$\widetilde{S}=\bigsqcup_{C \in \mathcal{C}'} C$$ for some $\mathcal{C}' \subset \mathcal{C}$ and $\mathcal{C}'=\psi^{-1}(T \setminus T')$ for some tableau $T' \subset T$, where $T=T(F,Q)$ and $\psi$ are as constructed in Definition~\ref{def: Q^2 tiling}. Now let $$\mathcal{E}= \left \{ \psi^{-1}(t_1,t_2) \text{ }| \text{ } (t_1,t_2) \in T\setminus T' \text{ and } (t_1-1,t_2) \notin T \setminus T' \text{ and } (t_1,t_2-1) \notin T \setminus T'   \right \}$$ denote the \textit{bottom-left} corners of $\mathcal{C}'$. Note that the element of $\mathcal{E}$ may be ordered vertically: we say $E$ is \textit{higher} than $E'$ if $t_2>t_2'$ where $E=\psi^{-1}(t_1,t_2)$ and $E'=\psi^{-1}(t'_1,t'_2)$. We now construct the $a_j$ recursively. Choose $a_1 \in E_1 \cap A$, where $E_1$ is the highest element of $\mathcal{E}$. Now suppose we have chosen $a_1, \ldots, a_j$ with each $a_j \in A\cap E_j$ for some $E_j \in \mathcal{E}$ (by minimality of $\widetilde{S}$, $A \cap E$ is non-empty for all $E \in \mathcal{E}$). We have that $y_r \leq a_j <y_{r+1}$ for some $r=r_j \in \{1, \ldots, |\mathcal{Y}|-1\}$ where $\mathcal{Y}=\mathcal{Y}(F,Q)=\{y_1<y_2< \ldots < y_{|\mathcal{Y}|} \}$ is as constructed in Definition~\ref{def: Q^2 tiling}. In fact $y_r$ and $y_{r+1}$ are the $y$ ordinates of the corners of the element of $\mathcal{E}$ that contains $a_j$. Now let $E_{j+1}$ be the highest element of $\mathcal{E}$ below the horizontal line $y= y_{r-1}$ and choose $a_{j+1} \in E_{j+1}$, and if such $E_{j+1}$ does not exist then $j=J$ and we are done with our construction. Now let $$G_1=(a_1+[0,\infty)^2) \cap F$$ and $$G_j=(a_j + [0,\infty)^2) \cap  F \cap \{(x,y) \in \mathbb{R} | y<y(a_{j-1})\}, $$ for $1<j\leq J$ . Note that $G_j$ is a tableau region translated by $a_j$. 

\begin{figure}[H]
\centering
\includegraphics[scale=0.4]{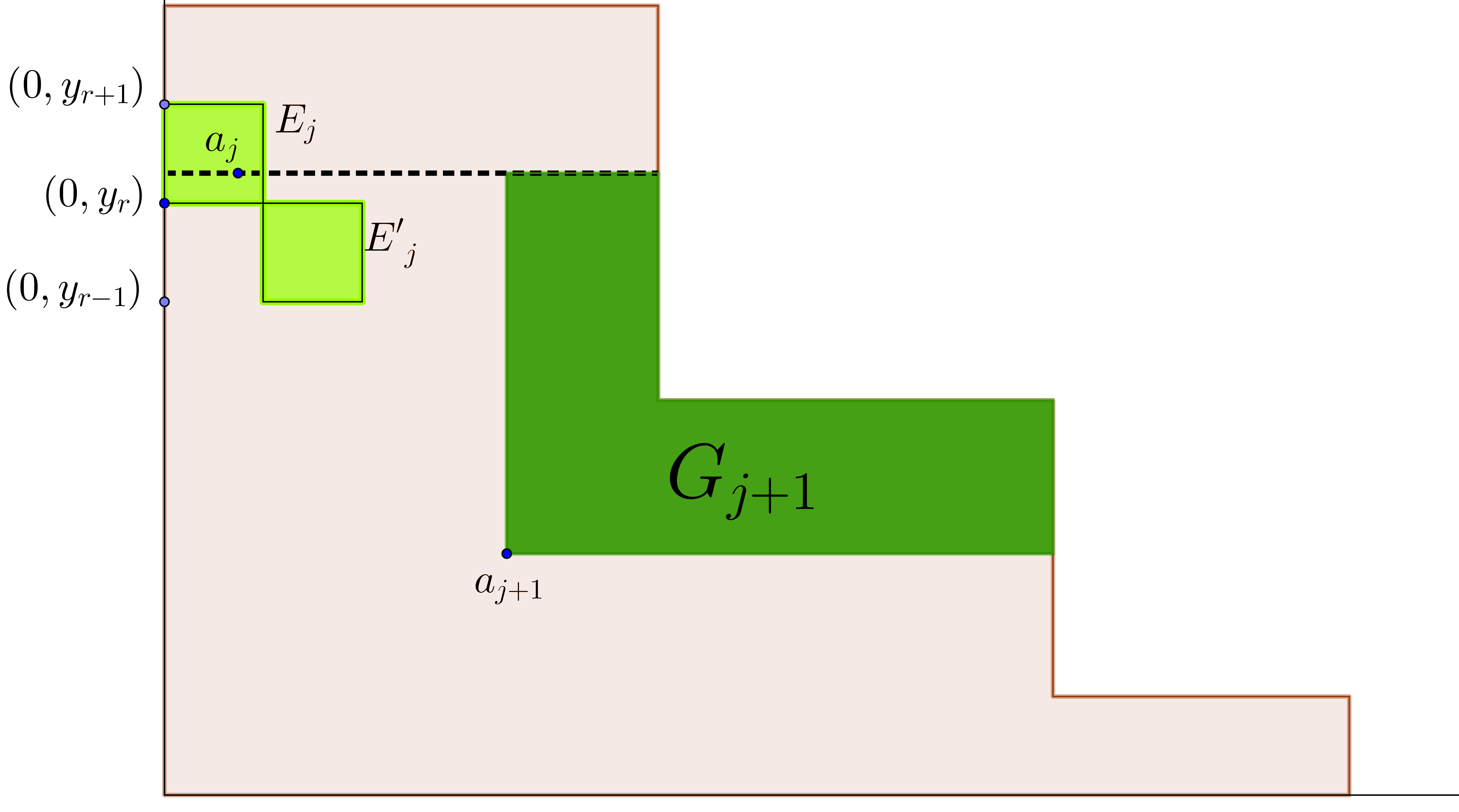}
\caption{ The tile $E'_j$ south-east of $E_j$ is an element of $\mathcal{E}$, but it is not $E_{j+1}$ since it is in the row of $\mathcal{C}$ just below $E_j$. This ensures the desired condition that $y(a_{j+1}) \leq y(a_j) - Q$. }
\label{fig: G_j construction}
\end{figure}

Now we let $$G=\bigsqcup_{j=1}^J G_j = \bigcup_{j=1}^J \left(a_j + [0,\infty)^2\right) \cap F.$$ As desired, we have that (cf. Figure~\ref{fig: G_j construction}) $$y(a_{j+1}) \leq y_{r-1} \leq y_r - Q \leq y(a_j) - Q.$$ 

\begin{figure}[H]
\centering
\includegraphics[scale=0.8]{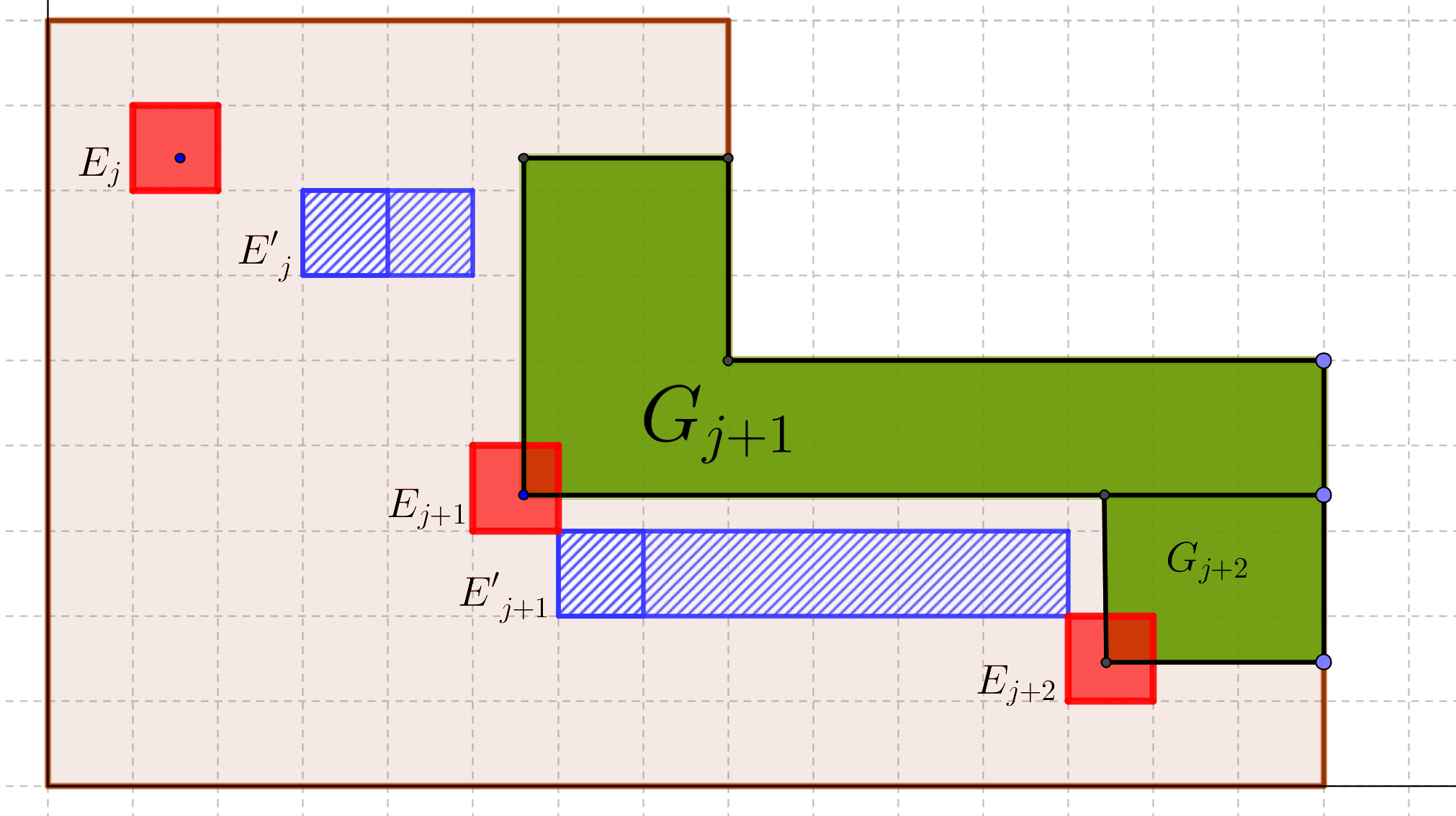}
\caption{ There is at most one element of $\mathcal{E}$ below $E_j$ and above $E_{j+1}$, which we call (assuming it exists) $E'_j$. The row of elements of $\mathcal{C}$ to the right of $E'_j$ that do not intersect $G$ is a subset of $\widetilde{S}_n \setminus \widetilde{G}$. The union of such rows is in fact exactly $\widetilde{S} \setminus \widetilde{G}$. These rows have disjoint projections onto the $x$-axis (one can see that the projection of $E_{j+1}$ onto the $x$-axis separates the projections of the two rows shown in the figure). Thus we get that the union of all such rows has Lebesgue measure at most $\frac{1}{Q} |F|$.}
\label{fig: tildeG estimate}
\end{figure}

It now remains to estimate $|S \setminus G|$. To this end, let $\widetilde{G}$ denote the smallest set that contains $G$ and may be written as a union of elements of $\mathcal{C}$. We have by Lemma~\ref{lemma: approximating regions} that \begin{align}\label{tildeG close to G} |\widetilde{G} \setminus G| \leq \frac{2|F|}{Q} \end{align} One may argue (see Figure~\ref{fig: tildeG estimate} and its caption below) that  \begin{align} \label{tildeS close to tildeG} |\widetilde{S} \setminus \widetilde{G}| \leq \frac{|F|}{Q}. \end{align}

Combining (\ref{tildeG close to G}) with (\ref{tildeS close to tildeG}) we get that $$ |S \setminus G| \leq |\widetilde{S} \setminus G| =|\widetilde{S} \setminus \widetilde{G}| + |\widetilde{G} \setminus G| \leq \frac{3|F|}{Q}$$ as desired. $\blacksquare$

We finish this section with a simple lemma which will allow us to remove a negligible set of integral points which lie too closely to the boundary of a tableaux region.

\begin{figure}[H]
\centering
\includegraphics[scale=0.2]{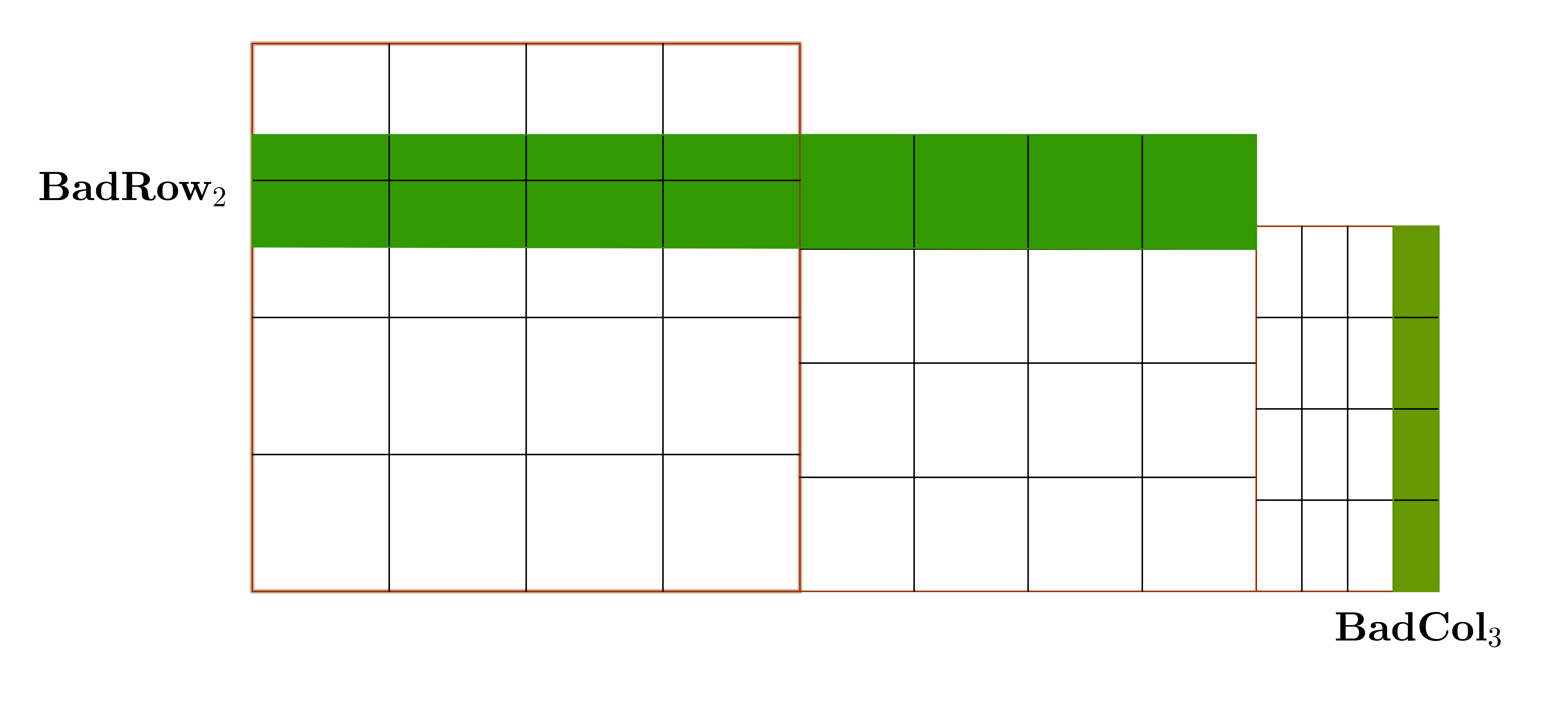}
\caption{An example of $\textbf{BadRow}_2$ and $\textbf{BadCol}_3$ with $Q=4$ and $\ell=3$. Note that $\textbf{BadRow}_2$ intersects exactly $\frac{1}{4}$ of $U_2$, strictly less than $\frac{1}{4}$ of $U_1$ and none of $U_3$.}
\label{fig: Bad Rows example}
\end{figure}

\begin{mydef}[Bad Rows and Bad Columns]\label{def: bad rows and columns} Fix the setup in Definition~\ref{def: Q^2 tiling} and partition $$F=\bigsqcup_{m=1}^{\ell} U_m$$ where $$U_m=[W_{m-1},W_m) \times [0,H_m)$$ for $m=1,2, \ldots, \ell$ with the convention $W_0=0$. Define (see Figure~\ref{fig: Bad Rows example}) for $m=1, \ldots, \ell$ the \textit{bad rows} $$\textbf{BadRow}_m=(U_1 \cup \ldots \cup U_m) \cap \left( \mathbb{R} \times \left[\frac{Q-1}{Q}H_m,H_m \right) \right)$$ and \textit{bad columns} $$\textbf{BadCol}_m= \bigsqcup_{j=0}^{Q-1} C^m_{Q,j} = \left [ W_{m-1} + \frac{Q-1}{Q}(W_m - W_{m-1}), W_{m}  \right) \times [0,H_m). $$

\end{mydef}

\begin{lemma}[Bad Rows and Columns removal]\label{bad row lemma} Fix the setup in Definition~\ref{def: bad rows and columns} and suppose $A' \subset F \cap \mathbb{N}^2$. Let $$A_0 = A' \setminus \left( \bigcup_{m=1}^{\ell} \textbf{BadRow}_m \cup \textbf{BadCol}_m \right).$$ Then $$|A_0| \geq |A'| - (\ell+1)Q^{-1}|F|.$$

\end{lemma}

\textbf{Proof:} We have that $|\textbf{BadRow}_m \cap U_j| \leq \frac{1}{Q}|U_j|$ for all $m,j \in \{1, \ldots, \ell\}$ (see Figure~\ref{fig: Bad Rows example}) and hence $|\textbf{BadRow}_m| \leq \frac{1}{Q}|F|$, which means that $$|\bigcup_{m=1}^{\ell} \textbf{BadRow}_m| \leq \frac{\ell}{Q}|F|.$$ On the other hand, $|\textbf{BadCol}_m|=\frac{1}{Q}|U_m|$. $\blacksquare$

\section{Proof of Theorem~\ref{main result in intro}}

\label{section: proof of main plunnecke inequality}

As in the hypothesis of Theorem~\ref{main result in intro}, fix $A,B \subset \mathbb{N}^2$ with $\alpha := \underline{d}_{\textbf{Tab}}(A)>0$ and $(0,0) \in B$ together with integers $0<k'<k$ and $L>0$. Now fix an integer \begin{align} \label{Q large} Q \geq 4(L+1) \alpha^{-1}\end{align} and choose (by Lemma~\ref{Lemma: sidelengths divide D}) a sequence $$F_n = \bigcup_{m=1}^{\ell(n)} [0,W_{n,m}) \times [0,H_{n,m})$$ with $\ell(n) \leq L$ such that 

\begin{enumerate}[label=(\roman*)] 
	\item $$\lim_{n \to \infty} \frac{|(A+kB) \cap F_n|}{|F_n|}  = \underline{d}_{\textbf{Tab}(L)}(A+kB)$$ 
	\item For all $n>1$ and $1 \leq m \leq \ell(n)$ we have $$H_{n,m} \equiv W_{n,m} \equiv 0 \mod Q^2$$ 	
	\item $$\lim_{n \to \infty} \inf_{m \in \{1, \ldots, \ell(n)\}} W_{n,m} = \lim_{n \to \infty} \inf_{m \in \{1, \ldots, \ell(m)\}} H_{n,m} = \infty.$$
\end{enumerate} 

We also assume, by reordering rectangles and deleting redundant rectangles if necessary, that $$H_{n,1}>H_{n,2}>\ldots >H_{n,\ell(n)}$$ $$W_{n,1}<W_{n,2}< \ldots < W_{n,\ell(n)}. $$ 

We now apply the techniques developed in Section~\ref{section: approximating subtableaux} to the tableau regions $F_n$. Note that the $\mathcal{C}(F_n,Q)$-measurable subtableau regions of $F_n$ are a union of at most $Q^2L^2$ (and thus a bounded function of $n$) rectangles with sidelengths tending to $\infty$ as $n \to \infty$ (see Remark~\ref{remark about sidelengths}), we thus have that \begin{align}\label{liminf alpha_n at least alpha} \liminf_{n \to \infty} \alpha_n \geq \alpha \end{align} where $$\alpha_n := \inf \left\{ \frac{|A \cap F'|}{|F'|} | F' \text{ is a } \mathcal{C}(F_n,Q) \text{-measurable subtableau region} \right\}.$$ Thus there exists $\tilde{N}$ (which depends on $Q$ and the sequence $F_n$, which we have fixed) such that $$\alpha_n > \frac{3}{4}\alpha \text{ for all } n>\tilde{N}.$$ From here on, we fix $n>\tilde{N}$ and let $\ell=\ell(n)$, $W_{m}=W_{n,m}$, $H_m=H_{n,m}$ for $m \in \{1, \ldots, \ell\}$. Applying Lemma~\ref{lemma: trimming set of points} to $F_n$ we obtain a subset $A' \subset A \cap F_n$ such that 
\begin{align}\label{trimmed above C-measurable tableau} \frac{|A' \cap (F_n \setminus F')|}{|F_n\setminus F'|} \leq \alpha_n + \frac{1}{Q^2} \text{ for all }\mathcal{C}(F_n,Q)\text{-measurable tableau } F' \subsetneqq F_n. \end{align} and $$\frac{|A'|}{|F_n|} \geq \alpha_n. $$

Now let $$A_0 = A' \setminus \left( \bigcup_{m=1}^{\ell} \textbf{BadRow}_m \cup \textbf{BadCol}_m \right)$$ where we have used the language of Definition~\ref{def: bad rows and columns}. We have by Lemma~\ref{bad row lemma} that \begin{align}\label{liminf of density of A_0} |A_0| \geq |A'| - (L+1)Q^{-1}|F_n| \geq (\alpha_n - (L+1)Q^{-1})|F_n|.\end{align} So from $\alpha_n> \frac{3}{4}\alpha$ and (\ref{Q large}) we get that \begin{align}\label{density of A_0 bounded away from zero} \frac{|A_0|}{|F_n|} > \frac{1}{2}\alpha \end{align} is bounded away from zero.

Applying the $\delta$-heavy truncated Pl\"unnecke's inequality (Proposition~\ref{truncated heavy plunnecke}) with $\delta=Q^{-1/2}$ to the finite set $A_0 \subset F_n$, we get \begin{align} \label{construction of A'_0} \frac{|(A +k'B) \cap F_n|}{|A_0|} \geq \frac{|(A_0+k'B) \cap F_n|}{|A_0|} \geq (1-Q^{-1/2})\left( \frac{|(A'_0+kB) \cap F_n|}{|A'_0|} \right) ^{k'/k} \end{align} for some nonempty $A'_0 \subset A_0$ with $$|A'_0| \geq Q^{-1/2}|A_0|.$$ Now let $$S_n = \bigcup_{v \in A'_0} (v+[0,\infty)^2) \cap F_n$$ be the smallest set that contains $A'_0$ and is the complement (in $F_n$) of a tableau region contained in $F_n$. We have by Lemma~\ref{lemma: approximating regions} that \begin{align} \frac{|\widetilde{S_n} \setminus S_n|}{|F|}\leq \frac{2}{Q} \end{align} where $\widetilde{S_n}$ is the smallest $\mathcal{C}$-measurable set that contains $S_n$.  Now note that \begin{align}\label{density of S_n bounded away from 0} \frac{|S_n|}{|F_n|} \geq \frac{|A'_0|}{|F_n|} \geq \frac{Q^{-1/2}|A_0|}{|F_n|} \geq \frac{1}{2}\alpha Q^{-1/2} \end{align} by (\ref{density of A_0 bounded away from zero}). Combining these two estimates gives \begin{align} \label{ratiowise tildeSn estimate} \frac{|\widetilde{S}_n|}{|S_n|} \leq 1+4\alpha^{-1} Q^{-1/2}. \end{align}

Also notice that $\widetilde{S}_n$ is the complement of a $\mathcal{C}(F_n,Q)$-measurable tableau, and thus we may apply (\ref{trimmed above C-measurable tableau}) to obtain 
\begin{align} \label{trimmed above tildeSn} \frac{|A'_0|}{|\widetilde{S}_n|} \leq \frac{|A' \cap \widetilde{S}_n|}{|\widetilde{S}_n|} \leq \alpha_n + \frac{1}{Q^2}. \end{align}

Now applying the inequalities (\ref{construction of A'_0}), followed by (\ref{liminf of density of A_0}), followed by (\ref{trimmed above tildeSn}) and then finally (\ref{ratiowise tildeSn estimate}) we obtain 

\begin{align*} 
& \left(  \frac{|(A+k'B) \cap F_n|}{|F_n|} \right)^{k/k'} \\ &\geq  \left(\frac{|A_0|}{|F_n|} (1-Q^{-1/2}) \right)^{k/k'} \frac{|(A'_0 +kB) \cap F_n|}{|A'_0|}
\\& \geq  \left((\alpha_n - (L+1)Q^{-1}) (1-Q^{-1/2}) \right)^{k/k'} \frac{|(A'_0 +kB) \cap F_n|}{|A'_0|}
\\& =  \left((\alpha_n - (L+1)Q^{-1}) (1-Q^{-1/2}) \right)^{k/k'} \frac{|\widetilde{S}_n|}{|A'_0|}\frac{|(A'_0 +kB) \cap F_n|}{|\widetilde{S}_n|} 
\\& \geq \left((\alpha_n - (L+1)Q^{-1}) (1-Q^{-1/2}) \right)^{k/k'} (\alpha_n + Q^{-2})^{-1}\frac{|(A'_0 +kB) \cap F_n|}{|\widetilde{S}_n|} 
\\& \geq \left((\alpha_n - (L+1)Q^{-1}) (1-Q^{-1/2}) \right)^{k/k'} (\alpha_n + Q^{-2})^{-1}(1+4\alpha^{-1}Q^{-1/2})^{-1} \frac{|(A'_0 +kB) \cap F_n|}{|S_n|} 
\end{align*}

which, by letting $\lambda(n,Q)$ denote the factor before $\frac{|(A'_0 +kB) \cap F_n|}{|S_n|} $, we rewrite as \begin{align} \label{def of lambda(n,Q)} \left(  \frac{|(A+k'B) \cap F_n|}{|F_n|} \right)^{k/k'} \geq \lambda(n,Q) \frac{|(A'_0 +kB) \cap F_n|}{|S_n|}. \end{align}

Using $\liminf_{n \to \infty} \alpha_n \geq \alpha$ (justified in (\ref{liminf alpha_n at least alpha})) it is an easy calculation to show that\footnote{We cannot reverse these limits, since the definiton of $\alpha_n$ depends on $Q$ and we also chose $n\geq \tilde{N}(Q)$. } \begin{align}\label{liminf of lambda(n,Q)} \liminf_{Q \to \infty} \liminf_{n \to \infty} \lambda(n,Q) \geq \alpha^{\frac{k}{k'} - 1}. \end{align} 

We now wish to show that $$ \liminf_{Q \to \infty} \liminf_{n \to \infty} \frac{|(A'_0 +kB) \cap F_n|}{|S_n|} \geq \underline{d}_{\textbf{Tab}(L)}(kB). $$

Note that this holds trivially in the case that $kB= \mathbb{N}^2$, as in this case we have that $|(A'_0 +kB) \cap F_n| = |S_n \cap \mathbb{Z}^2|$ and thus \begin{align}\label{special case when kB=N^2} \frac{|(A'_0 +kB) \cap F_n|}{|S_n|}  = 1. \end{align}  We now return to the general case. Apply Lemma~\ref{Lemma: constructing a_j} to $A'_0$ and $S_n$ to obtain $a_1,a_2, \ldots, a_J \in A$ with\footnote{Recall the notation $y(a,b)=b$.} $$y(a_j) \leq y(a_{j-1}) - Q$$ such that the set $$G:= \bigcup_{j=1}^J \left(a_j + [0,\infty)^2\right) \cap F_n$$ satisfies \begin{align}\label{G close to S} \frac{|S \setminus G|}{|F_n|} \leq \frac{3}{Q}.\end{align}

Now decompose $$G= \bigsqcup_{j=1}^J G_j$$ where $$G_1=(a_1+[0,\infty)^2) \cap F_n$$ and $$G_j=(a_j + [0,\infty)^2) \cap  F_n \cap \{(x,y) \in \mathbb{R} | y<y(a_{j-1})\}, $$ for $1<j\leq J$. Note that $G_j$ is a tableau region translated by $a_j$ (cf. Figure~\ref{fig: G_j construction} in the proof of Lemma~\ref{Lemma: constructing a_j}). 

\textbf{Claim:} Each $G_j$ is a union of at most $L$ rectangles with bottom corner $a_j$, each with sidelengths at least $Q$.

\textbf{Proof of Claim:} Writing $a_j = (x,y)$, we can write (see Figure~\ref{fig: G_j construction}) $$G_j= \bigcup_{m} [x,W_m) \times [y, \min\{ y(a_{j-1}), H_m\} )  $$ where the union is over $m \in \{1, \ldots, \ell\}$ such that $x <W_m$ and $y<H_m$ (for $j=1$, we omit the $y(a_{j-1})$). For such $m$, we have that $W_m-x \geq Q$ and $H_m-y \geq Q$ since $a_j$ avoids $\textbf{BadRow}_m$ and $\textbf{BadCol}_m$, respectively. We also have $y(a_{j-1})-y \geq Q$ by construction of the $a_j$. This completes the proof of the Claim. $\square$

This claim implies that $$\frac{|(A'_0+kB) \cap G_j|}{|G_j|}\geq \frac{|(a_j+kB) \cap G_j|}{|G_j|} \geq \underline{d}_{\textbf{Tab}(L)}(kB) -\delta(Q)$$ where $$\delta:\mathbb{N} \to [0,\infty)$$ is a map (it depends on $kB$, but not $n$) such that\footnote{To see this, apply Lemma~\ref{characterization of lower tableaux density} to $kB$ with $R=Q$ and $\epsilon=\delta(Q)$.} $\lim_{Q \to \infty} \delta(Q)=0$. We thus have that $$\frac{|(A'_0 +kB) \cap G|}{|G|} \geq \underline{d}_{\textbf{Tab}(L)}(kB) -\delta(Q).$$ Now combining (\ref{G close to S}) with (\ref{density of S_n bounded away from 0}) gives $$\frac{|G|}{|S_n|} \geq 1 - 6\alpha^{-1}Q^{-1/2}$$ and so we have that \begin{align*} \frac{|(A'_0 + kB) \cap S_n|}{|S_n|} \geq & \frac{|(A'_0 + kB) \cap G|}{|S_n|} \\ \geq & \frac{|G|}{|S_n|} (\underline{d}_{\textbf{Tab}(L)}(kB) -\delta(Q)) \\ \geq & \left( 1 - 6\alpha^{-1}Q^{-1/2} \right)(\underline{d}_{\textbf{Tab}(L)}(kB) -\delta(Q)).  \end{align*} So in summary, we have (see (\ref{def of lambda(n,Q)})) that $$\left(  \frac{|(A+k'B) \cap F_n|}{|F_n|} \right)^{k/k'} \geq \lambda(n,Q)  \left( 1 - 6\alpha^{-1}Q^{-1/2} \right)(\underline{d}_{\textbf{Tab}(L)}(kB) -\delta(Q)).$$ Finally, this completes the proof as \begin{align*} \underline{d}_{\textbf{Tab}(L)}(A+k'B) & = \liminf_{n \to \infty} \frac{|(A+k'B) \cap F_n|}{|F_n|} \\ & \geq \liminf_{n \to \infty} \lambda(n,Q)^{k'/k} \left(1 - 6\alpha^{-1}Q^{-1/2} \right)^{k'/k} (\underline{d}_{\textbf{Tab}(L)}(kB) -\delta(Q))^{k'/k} \end{align*} from which we deduce Theorem~\ref{main result in intro} by letting $Q \to \infty$ and using (\ref{liminf of lambda(n,Q)}). $\blacksquare$

\section{Further questions}

\label{sec: questions}

We now list some related open problems. We start by recalling the main motivating question of this paper.

\begin{question} For $A,B \subset \bN^2$ and positive integers $k'<k$, is it true that $$\underline{d}_{\textbf{Rect}}(A+k'B) \geq \underline{d}_{\textbf{Rect}}(A)^{1 - \frac{k'}{k}}\underline{d}_{\textbf{Rect}}(kB)^{\frac{k'}{k}}?$$

\end{question}

One may state a finitistic version of this problem by considering the finitistic and multidimensional analogue of the classical Schnirelmann density given by $$\sigma_{N,M}(A)=\min_{0\leq n \leq N, 0 \leq m \leq M} \frac{A \cap ([0,n] \times [0,m])}{(n+1)(m+1)}$$ for $A \subset \bN^2 \cap ([0,N] \times [0,M])$.

\begin{question}[Multidimensional Schnirelmann density Pl\"unnecke inequality]\label{question: schnirelmann multidim} Fix positive integers $N$ and $M$. If $A,B \subset \bN^2 \cap ([0,N] \times [0,M])$, with $(0,0) \in B$, then is it true that $$ \sigma_{N,M}(A+k'B) \geq \sigma_{N,M}(A)^{1 - \frac{k'}{k}}\sigma_{N,M}(kB)^{\frac{k'}{k}}$$ for positive integers $k'<k$ ?

\end{question}

We note that the $M=0$ case is precisely Pl\"unnecke's classical inequality for Schnirelmann density in $\bN$ (Theorem~\ref{Plunnecke for Schnirelmann}). Even the following special case is not clear. 

\begin{question} What is the answer to Question~\ref{question: schnirelmann multidim} in the case $M=1$?

\end{question}

\appendix

\section{Densities of fractal sets}

\label{section: Densities of fractal sets}

We will now prove our formulae for $\underline{d}_{\textbf{Tab}}(A)$ and $\underline{d}_{\textbf{Rect}}(A)$ (Proposition~\ref{formula for density of fractal}) for fractal sets $A$. We will use the notation $[a,b)=\{ x \in \mathbb{Z} | a \leq x < b\}$. So let us fix in this section $A=A(P,N)$ a fractal set, where $N \in \mathbb{N}$ and $P \subset \{0,1, \ldots, N\}^2$ is a pattern. We also fix the data $A_k \subset \cP_k$ and $u_k \in \mathbb{N}$ given in Definition~\ref{def: fractal set}. 

\subsection{Rectangular density}

\label{subsection: rectangular density}

We first deal with the rectangular case; so in this subsection we fix a F{\o}lner sequence of the form $F_n=[0,W_n) \times [0,H_n) \cap \mathbb{N}^2$, $n=1,2,3 \ldots$ where $W_n,H_n \to \infty$. It is enough to show that \begin{align}\label{appendix: desired inequality for rectangular} \liminf_{n \to \infty} \frac{|A \cap F_n|}{|F_n|} \geq \alpha :=  \min \left\{ \frac{|P \cap I|}{|I|} \text{ } | \text{ } I=[0,x] \times [0,y] \cap \mathbb{N}^2 \text{ for some } x,y \in \{0,1, \ldots, N\}  \right \} \end{align} as it is easy to see that $$\alpha \geq \underline{d}_{\textbf{Tab}}(A)$$ by considering the F{\o}lner sequence $$F_k=[0,(x+1)u_k) \times [0,(y+1)u_k)$$ where $(x,y) \in \{0,1 \ldots,N\}^2$ attain the minimum in (\ref{appendix: desired inequality for rectangular}).

Define the \textit{core} of $F_n$ to be the set $$F^o_n = F_n \cap [0,(N+1)u_{k_n})^2 $$ where $k_n$ is the largest positive integer such that $[0,u_{k_n})^2 \subset F_n$. Our goal is to show the following two inequalities 
\begin{align} \label{F^o asymptotic} \liminf_{n \to \infty} \frac{|F^o_n \cap A|}{|F^o_n|} \geq \alpha \end{align}

\begin{align} \label{F^o complement asymptotic} \liminf_{n \to \infty} \frac{|(F_n\setminus F^o_n) \cap A|}{|F_n \setminus F^o_n|} \geq \alpha \end{align}

which will imply the formula for rectangular density. The key ingredient is the following \textit{Perturbation Lemma}.

\begin{lemma}[Perturbation Lemma] \label{Lemma: perturbation}Suppose $F=[0,W) \times [0, H) \cap \mathbb{N}^2$ where $(W,H) \in \mathbb{N}^2$. Suppose that for some $k \in \mathbb{N}$ we have that $u_k \leq H \leq u_{k+1}$ and $W \leq u_{k+1}$. Then there exists $\tilde{H} \in \{ iu_k | i \in \{1,2, \ldots N+1 \} \}$ such that $$\frac{|\tilde{F} \cap A|}{|\tilde{F}|} \leq \frac{|F \cap A|}{|F|}$$ where $$\tilde{F}=[0,W) \times [0,\tilde{H}) \cap \mathbb{N}^2.$$

\end{lemma}

\begin{figure}[H]
\includegraphics[scale=0.56]{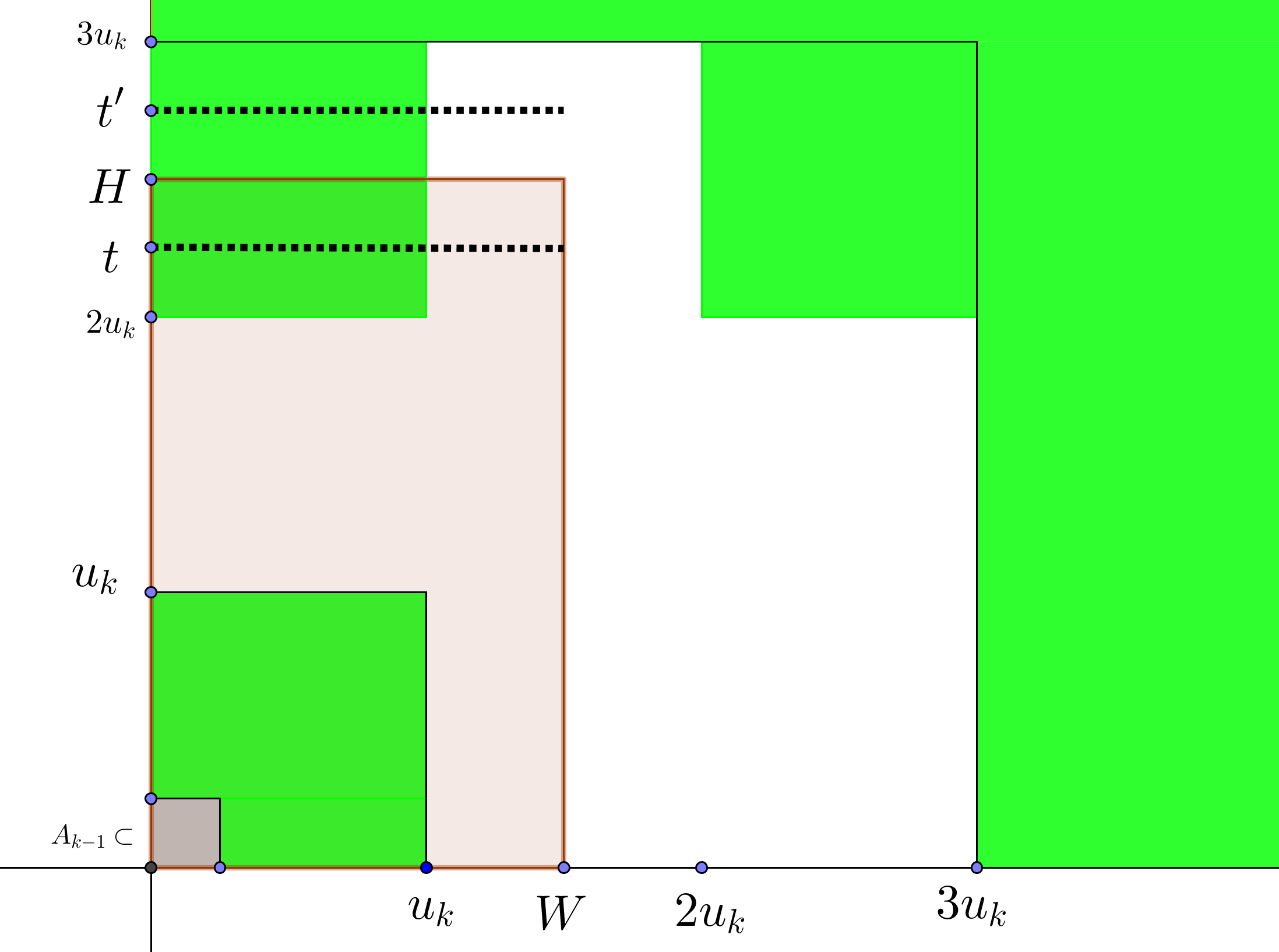}
\centering
\caption{The cross sections $A \cap ([0,W) \times \{t \})$ and $A \cap ([0,W) \times \{t' \})$ are the same upto a vertical translation.}
\label{Fig: cross sections}

\end{figure}

\textbf{Proof:} Since $F \setminus [0,(N+1)u_k)^2 \subset A$, it suffices to consider the case where $W,H \leq (N+1)u_k$. Write $H=q u_k + r$ where $r \in \{0,1, \ldots u_k-1\}$ and $q \in \{1, \ldots, (N+1)\}$. If $q=N+1$ then we must have $r=0$, which means that we can set $\tilde{H}:=H=(N+1)u_k$ and we are done with the proof, so let us now assume $q<N+1$. The quantity $$ \frac{|A \cap ([0,W) \times \{t \})|}{|W|} $$ is constant for $t \in \{0,1, \ldots u_k - 1\}$ since the sets $A \cap \left( [0,W) \times \{t \} \right)$ all have the same projection onto the $x$-axis (see Figure~\ref{Fig: cross sections}). We denote this constant by $C$. 

It will be convenient to use the notation $$d(X|Y)=\frac{|X \cap Y|}{|Y|}$$ where $X,Y$ are finite sets. Hence if we define, for $t=0,1, \ldots u_k$, the expression $$g(t) := d\left( A | [0,W) \times [0,qu_k +t) \right) = \frac{|A \cap [0,W) \times [0,qu_k + t)|}{W(qu_k+t)}$$ then we can rewrite it as a convex combination \begin{align} g(t) = \lambda \cdot d\left(A| [0,W) \times [0,qu_k) \right) +(1-\lambda)C  \end{align} where $$\lambda = \frac{Wqu_k}{W(qu_k+t)}. $$  Now notice that $g(t)$ is a monotonic function of $\lambda$, and thus a monotonic function of $t$. Hence the minimum of $g(t)$ occurs at either $t=0$ or $t=u_k$. If it occurs at $t=0$, we set $\tilde{H}=qu_k$ while if it occurs at $t=u_k$, we set $\tilde{H}=(q+1)u_k$. We thus get that $$\frac{|A \cap \tilde{F}|}{|\tilde{F}|} \leq \frac{|F \cap A|}{|F|} $$ where $\tilde{F}=[0,W) \times [0,\tilde{H}).$ $\blacksquare$

\begin{remark}\label{remark: perturbation} In Lemma~\ref{Lemma: perturbation} one may drop the hypothesis that $u_{k+1} \geq H$ if one replaces $A$ with $\cP_k$. The reason is that this hypothesis was only used to show that the cross section density $C$ was constant with respect to $t$, which is not an issue if one replaces $A$ with $\cP_k$. This variation of the lemma will also be useful.

\end{remark}

\begin{lemma}[Perturbing both sides] \label{lemma: perturbing both sides} In the situation of Lemma~\ref{Lemma: perturbation} if it is also the case that $W \geq u_k$, then there exists $\tilde{H}, \tilde{W} \in \{ iu_k \text{ }|\text{ } i \in \{1,2, \ldots N+1 \} \}$ such that $$\frac{|\tilde{F} \cap A|}{|\tilde{F}|} \leq \frac{|F \cap A|}{|F|}$$ where $$\tilde{F}=[0,\tilde{W)} \times [0,\tilde{H}) \cap \mathbb{N}^2.$$

\end{lemma}

\textbf{Proof:} After applying Lemma~\ref{Lemma: perturbation}, apply it again with dimensions reversed. $\blacksquare$

We may apply Lemma~\ref{lemma: perturbing both sides} with $F=F^o_n$ (and $k=k_n$) to obtain the corresponding rectangle $\tilde{F}^o_n$ which satisfies $$\frac{|\tilde{F}^o_n \cap A|}{|\tilde{F}^o_n|} \leq \frac{|F^o_n \cap A|}{|F^o_n|}.$$ However, by definition of $\alpha$, it is clear that $$\frac{|\tilde{F}^o_n \cap \cP_{k_n}|}{|\tilde{F}^o_n|} \geq \alpha.$$ But since $\frac{u_k}{u_{k-1}} \to \infty$, we have that $$ |\frac{|\tilde{F}^o_n \cap \cP_{k_n}|}{|\tilde{F}^o_n|} - \frac{|\tilde{F}^o_n \cap A|}{|\tilde{F}^o_n|} | \to 0 \text{ as } n \rightarrow \infty,$$ which shows (\ref{F^o asymptotic}). Now we turn to showing (\ref{F^o complement asymptotic}). Suppose, WLOG, that $H_n \geq W_n$, thus $$F_n \setminus F^o_n = [0,W_n) \times [(n+1)u_{k_n},H_n).$$ Now if we are in the degenerate case where $H_n \leq u_{k_n+1}$ then in fact $F_n \setminus F^o_n= \emptyset$ and so we are done. If however $H_n > u_{k_n+1}$ then $\frac{|F_n \setminus F^o_n|}{|F_n|} \geq 1 - \epsilon_n$, for some $\epsilon_n \to 0$, which means that we may replace $F_n \setminus F^o_n$ with $F_n$ in (\ref{F^o complement asymptotic}) to get a logically equivalent statement. So we apply Lemma~\ref{Lemma: perturbation} (with $k$ being the largest integer such that $H_n \geq u_{k}$ and $F=F_n$) and argue as before.

\begin{remark} We did not really need to consider the definition of a \textit{core}, but this notion will be useful in establishing the formula for Tableaux density in the next subsection.

\end{remark}

\subsection{Tableaux density}

We will now prove the formula for $\underline{d}_{\textbf{Tab}}(A)$ for the fractal set $A=A(N,P)$. In this subsection we fix an integer $L \geq 1$ and a F{\o}lner sequence $$ F_n = \bigcup_{j=1}^L \left [0,W_{n,j}\right) \times \left[0,H_{n,j}\right). $$ Now let $$\beta=\min\left \{ \frac{|P \cap T|}{|T|} \text{ }| \text{ } T \subset \{0,1, \ldots, N\}^2 \text{ is a tableaux } \right\}.$$ We wish to show that $$\beta \leq \liminf_{n \to \infty} \frac{|A \cap F_n|}{|F_n|}, $$ where $H_{n,j}$ is decreasing and $W_{n,j}$ is increasing in $j$, for each fixed $n$.

Exactly as in Section~\ref{subsection: rectangular density}, we define the \textit{core} of $F_n$ to be the set $$F^o_n = F_n \cap [0,u_{1+k_n})^2 $$ where $k_n$ is the largest positive integer such that $[0,u_{k_n})^2 \subset F_n$. As before, we start with a perturbation lemma. We note again that (since $\frac{u_k}{u_{k-1}}\to \infty$) \begin{align} \lim_{n \to \infty} \left| \frac{|A \cap F_n|}{ |F_n|} - \frac{|\cP_{k_n} \cap F_n|}{|F_n|}\right| = 0. \end{align}

\begin{lemma}[Order $L$ perturbation lemma] \label{lemma: tab perturbation} Let $$F= \bigcup_{j=1}^L [0,W_j) \times [0,H_j)$$ for some positive integers $H_1 \geq \cdots \geq H_L$ and $W_1 \leq \cdots W_L$. Now suppose that there exists a positive integer such that $k$ be the largest positive integer such that  $[0,u_k)^2 \subset F \subset [0,u_{k+1})^2$. Then, for each $j=1, \ldots, L$, there exists $\tilde{W}_j, \tilde{H}_j \in \{iu_{k} \text{ }|\text{ } i \in \{1,\ldots,N+1\}\}$ such that $$\frac{|\cP_k \cap \tilde{F}|}{|\tilde{F}|} \leq \frac{|\cP_k \cap F|}{|F|}$$ where $$\tilde{F}=\bigcup_{j=1}^L [0,\tilde{W}_j) \times [0,\tilde{H}_j)$$ and $\tilde{H}_1 \geq \cdots \geq \tilde{H}_L$ and $\tilde{W}_1 \leq \cdots \leq \tilde{W}_L$

\end{lemma}

\textbf{Proof (Sketch):} The $L=1$ case is essentially Lemma~\ref{lemma: perturbing both sides} (see Remark~\ref{remark: perturbation}. In fact, in this proof we always use the formulation given in that remark). For $L>1$ we proceed by induction as follows. Let $H=H_1,W=W_1$ and apply Lemma~\ref{Lemma: perturbation} to $[0,H) \times [0,W)$ to produce $\tilde{H}$. Now if $\tilde{H}\geq H_2$ then set $\tilde{H}_1=\tilde{H}$. However if $\tilde{H}_1< H_2$ then we see, by the same convexity argument as in Lemma~\ref{Lemma: perturbation}, that $$\frac{|\cP_k \cap F'|}{|F'|} \leq \frac{|\cP_k \cap F|}{|F|}$$ where $$F'=\bigcup_{i=1}^{L-1} [0,W_{i+1}) \times [0,H_{i+1})$$ and so we are done (in this case) by applying the induction hypothesis to $F'$. We construct $\tilde{H}_2, \ldots \tilde{H}_L$ by continuing in this way (for example, to get $\tilde{H}_2$ one applies the same technique to $[W_1,W_2) \times [0,H_2)$).   $\blacksquare$

The formula for the tableaux density of $A$ now follows by applying Lemma~\ref{lemma: tab perturbation} to the core of $F_n$ (with $k=k_n$) and arguing as we did in the case of rectangular density. 

\section{Density of cartesian products}

\label{appendix: density of cartesian products}

We now prove the property \begin{align}\label{appendix: product formula density} \underline{d}_{\textbf{Tab}(L)} (A \times B) = \underline{d}(A) \underline{d}(B) \end{align} for $A,B \subset \bN$ and $L \in \bZ_{>0}$. It is clear that  $$\underline{d}_{\textbf{Tab}(L)} (A \times B) \leq \underline{d}(A) \underline{d}(B)$$ and so we focus on proving the reverse inequality.  We do this by induction on $L$. The $L=1$ case is clear, so let us suppose that $L>0$ and that the property holds for $\textbf{Tab}(L-1)$. Let $A,B \subset \bN$ and let $$F_n = \bN^2 \cap \bigcup_{i=1}^L [0,W_{n,i}) \times [0,H_{n,i})$$ be a sequence (where $W_{n,1} < \cdots < W_{n,L}$ and $H_{n,1} > \cdots > H_{n,L}$ are positive integers with $W_{n,i}, H_{n,i} \to \infty$ as $n \to \infty$ for each $i$) such that $$\underline{d}_{\textbf{Tab}(L)}(A \times B) = \lim_{n \to \infty} \frac{|(A \times B) \cap F_n|}{|F_n|}.$$  Now we assume, by passing to a subsequence if necessary, that the limit $$\alpha_i := \lim_{n \to \infty} \frac{|A \cap [W_{n,i},W_{n,i+1})|}{W_{n,i+1}-W_{n,i}}$$ exists for $i=0, \ldots, L-1$ (where $W_{n,0}=0$). Likewise, we assume that $$\beta_j:= \lim_{n \to \infty} \frac{|B \cap [H_{n,j+1},H_{n,j+2})|}{H_{n,j+2}-H_{n,j+1}}$$ exists for $j=0, \ldots, L-1$ (where $H_{n,L+1}=0$). Using the fact that $\alpha_0 \geq \underline{d}(A)$ and $\beta_{L-1} \geq \underline{d}(B)$ we can easily deduce that one of the following must occur:

\begin{enumerate}[label=(\textbf{\alph*})]
	\item \label{case remove} There exists $r \in \{0, \ldots ,L-1\}$ such that $\alpha_r \geq \underline{d}(A)$ and $\beta_r \geq \underline{d}(B)$.
	\item \label{case add}There exists $r \in \{0, \ldots, L-2\}$ such that $\alpha_{r+1} < \underline{d}(A)$ and $\beta_r < \underline{d}(B)$.

\end{enumerate}

Let us first consider the case where \ref{case remove} occurs. Hence \begin{align}\label{edge box is dense} \lim_{n \to \infty} \frac{|(A \times B) \cap \left([W_{n,r}, W_{n,r+1}) \times [H_{n,r+2},H_{n,r+1})\right)|}{|[W_{n,r}, W_{n,r+1}) \times [H_{n,r+2},H_{n,r+1})|} = \alpha_r\beta_r\geq \underline{d}(A)\underline{d}(B). \end{align} Now consider the sequence $$F'_n=F_n \setminus \left( [W_{n,r}, W_{n,r+1}) \times [H_{n,r+2},H_{n,r+1}) \right) = \bN^2 \cap \bigcup_{i \neq r} [0,W_{n,i}) \times [0,H_{n,i})$$ and observe that it is in $\textbf{Tab}(L-1)$. Thus the induction hypothesis implies that \begin{align} \label{invoking induction L-1} \liminf_{n \to \infty} \frac{|(A \times B) \cap F'_n|}{|F'_n|} \geq \underline{d}_{\textbf{Tab}(L-1)} = \underline{d}(A)\underline{d}(B).\end{align} Combining (\ref{edge box is dense}) and (\ref{invoking induction L-1}) give $$ \liminf_{n \to \infty} \frac{|(A \times B) \cap F_n|}{|F_n|} \geq \underline{d}(A) \underline{d}(B), $$ which completes the induction step in this case.

Now suppose that \ref{case add} occurs. Hence \begin{align} \label{skinny case b} \lim_{n \to \infty} \frac{|(A \times B) \cap \left([W_{n,r+1}, W_{n,r+2}) \times [H_{n,r+2},H_{n,r+1})\right)|}{|[W_{n,r+1}, W_{n,r+2}) \times [H_{n,r+2},H_{n,r+1})|} = \alpha_{r+1}\beta_r<\underline{d}(A)\underline{d}(B). \end{align}

Now, let $$F''_n = F_n \sqcup \left([W_{n,r+1}, W_{n,r+2}) \times [H_{n,r+2},H_{n,r+1})\right) \cap \bN^2,$$ and observe that the sequence $\left(F''_n\right)_n$ is an element of $\textbf{Tab}(L-1)$. Hence the induction hypothesis implies that \begin{align} \label{induction case b} \liminf_{n \to \infty} \frac{|(A \times B) \cap F''_n|}{|F''_n|} \geq \underline{d}(A) \underline{d}(B).\end{align} Combining the estimates (\ref{skinny case b}) and (\ref{induction case b}) gives $$ \liminf_{n \to \infty} \frac{|(A \times B) \cap F_n|}{|F_n|} \geq \underline{d}(A) \underline{d}(B), $$ which completes the induction step in this case.

\bibliographystyle{plain}
\bibliography{mybib}

\begin{thebibliography}{1}

\bibitem{BjorklundFishPlunnecke}
M.~Bj{\"o}rklund and A.~Fish.
\newblock Pl{\"u}nnecke inequalities for countable abelian groups.
\newblock {\em Accepted in Journal fur die Reine und Angewandte Mathematik
  (Crelle)}, 2013.

\bibitem{BulinskiFishPlunnecke}
K.~Bulinski and A.~Fish.
\newblock Pl\"unnecke inequalities for measure graphs with applications.
\newblock {\em Accepted in Ergodic Theory and Dynamical Systems, preprint:
  arXiv:1407.4174v2}, 2014.

\bibitem{ErdosBasis}
P.~Erd{\H{o}}s.
\newblock On the arithmetical density of the sum of two sequences one of which
  forms a basis for the integers.
\newblock {\em Acta Arithmetica}, 1(2):197--200, 1935.

\bibitem{JinPlunnecke}
R.~Jin.
\newblock Pl\"unnecke's theorem for asymptotic densities.
\newblock {\em Trans. Amer. Math. Soc.}, 363(10):5059--5070, 2011.

\bibitem{JinEpsilon}
R.~Jin.
\newblock Density versions of {P}l\"unnecke inequality: epsilon-delta approach.
\newblock In {\em Combinatorial and additive number theory---{CANT} 2011 and
  2012}, volume 101 of {\em Springer Proc. Math. Stat.}, pages 99--113.
  Springer, New York, 2014.

\bibitem{Plunnecke1970}
H.~Pl{\"u}nnecke.
\newblock Eine zahlentheoretische {A}nwendung der {G}raphentheorie.
\newblock {\em J. Reine Angew. Math.}, 243:171--183, 1970.

\bibitem{Ruzsasumsetsandstructure}
I.~Z. Ruzsa.
\newblock Sumsets and structure.
\newblock In {\em Combinatorial number theory and additive group theory}, Adv.
  Courses Math. CRM Barcelona, pages 87--210. Birkh\"auser Verlag, Basel, 2009.

\bibitem{SchnirelmannEigenschaften}
L.~Schnirelmann.
\newblock \"{U}ber additive {E}igenschaften von {Z}ahlen.
\newblock {\em Math. Ann.}, 107(1):649--690, 1933.

\end{thebibliography}

\end{document}